\documentclass[11pt]{article}
\usepackage[centertags]{amsmath}
\usepackage{amsmath}
\usepackage{amssymb}
\usepackage{amsthm}
\numberwithin{equation}{section}

\newcommand{\be}{\begin{equation}}
\newcommand{\ee}{\end{equation}}

\newtheorem{definition}{Definition}[section]
\newtheorem{theorem}{Theorem}[section]
\newtheorem{proposition}{Proposition}[section]
\newtheorem{lemma}{Lemma}[section]
\newtheorem{remark}{Remark}[section]
\newtheorem{corollary}{Corollary}[section]

\newcommand{\ind}{\mathbf{1}}

\newcommand{\eps}{\varepsilon}

\def \esssup {\mbox{ess sup}}
\def \nn{\nonumber}
\def \cp {{\cal P}}

\def \R{\mathbb{R}}
\def \P{\mathbb{P}}

\def \E{\mathbb{E}}

\def \H{\mathbb{H}}

\def \bf{\textbf}
\def \it{\textit}
\def \sc{\textsc}

\def \nd {\noindent}

\def \F{\mathcal{F}}
\def \g{\mathcal{g}}

\def \lb {\label}
\def \cs {{\cal S}}
\def \mx{\mbox}
\def \G{\Gamma}
\def \cu{{\cal U}}
\def \ba{\begin{array}}
\def \ea {\end{array}}
\def \bal{\begin{array}{l}}
\def \ball{\begin{array}{ll}}
\def \a {\alpha}

\usepackage{amssymb}
\usepackage{amsthm}
\usepackage{amsmath}
\usepackage{latexsym}
\usepackage[T1]{fontenc}
\usepackage[latin1]{inputenc}

\def \gm {\gamma}

\def \R{\mathbb{R}}
\def \P{\mathbb{P}}

\def \E{\mathbb{E}}

\def \H{\mathbb{H}}

\def \bf{\textbf}
\def \it{\textit}
\def \sc{\textsc}

\def \ni {\noindent}

\def \F{\mathcal{F}}
\def \g{\mathcal{g}}

\numberwithin{equation}{section}
\def \cl {\mathcal{L}}

\def \R{\mathbb{R}}
\def \P{\mathbb{P}}

\def \E{\mathbb{E}}

\def \H{\mathbb{H}}

\def \bf{\textbf}
\def \it{\textit}
\def \sc{\textsc}

\def \ni {\noindent}

\def \F{\mathcal{F}}
\def \g{\mathcal{g}}

\def \o {\omega}
\def \O {\Omega}
\def \ms {\medskip}
\def \bs {\bigskip}
\def \ed {\end{document}}
\def \cc {{\cal C}}
\def \ca {{\cal A}}
\def \fr {\forall}
\def \t {\tau}
\def \frt {\forall t\le T}
\def \txt {(t,x)\in [0,T]\times \R^d}
\def \esp {[0,T]\times \R^d}
\def \tx {{t,x}}
\def \txp{(t,x)}
\def \cu {{\cal U}}
\def \rw {\rightarrow}

\def \g{\gamma}
\def \l{\lambda}
\def \spo{[0,T)\times \R^d}
\def \btheta {\bar \theta}
\def\eps {\epsilon}
\def \lb{\label}

\def \cC{\cc}
\def \vro{\varrho}
\def \pg {\Pi_{\mbox{pg}}}

\begin{document}
\title{On the Stochastic Control-Stopping Problem }
\author{ Brahim EL ASRI \thanks{Universit\'e Ibn Zohr, Equipe. Aide \`a la decision,
ENSA, B.P.  1136, Agadir, Maroc. e-mail: b.elasri@uiz.ac.ma } \,,\,\,Said HAMADENE \thanks{LMM, Le Mans Universit\'e, Avenue Olivier Messiaen, 72085 Le Mans, Cedex 9, France. e-mail:
hamadene@univ-lemans.fr}\,\,\,
\, and \, Khalid OUFDIL \thanks{Universit\'e Ibn Zohr, Equipe. Aide \`a la decision,
ENSA, B.P.  1136, Agadir, Maroc. e-mail: khalid.ofdil@gmail.com.}}\maketitle

\begin{abstract}
We study the stochastic control-stopping problem when the data are of polynomial growth. The approach is based on 
backward stochastic differential equations (BSDEs for short). The problem turns into the study of a specific reflected BSDE with a stochastic Lipschitz coefficient for which we show existence and uniqueness of the solution. We then establish its relationship with the value function of the control-stopping problem. The optimal strategy is exhibited. Finally in the Markovian framework we prove that the value function is the unique viscosity solution of the associated Hamilton-Jacobi-Bellman equation.
\end{abstract}
\noindent {${Keywords}:$} Reflected Backward stochastic differential
equations; Mixed stochastic control; control-stopping problem; stochastic Lipschitz condition;  Hamilton-Jacobi-Bellman equation; Viscosity solution.
\ms

\section{Introduction}
The main objective of this paper is to deal with the finite horizon  control-stopping problem in its weak formulation (see e.g. \cite{B,DE, dava, hlp00} to quote a few, and the references therein) which we describe hereafter. 

Let us consider a controlled system whose dynamics $(x_s)_{s\le T}$ is a weak solution of the following functional stochastic differential equation: 
\begin{equation}
dx_s=f(s,x,u_s)ds+\sigma(s,x)dB^u_s, \,s\in[0,T] \mx{ and } x_0\in \R^d \mx{ fixed};
\end{equation}
$u:=(u_s)_{s\le T}$ is a stochastic process by which the controller intervenes on the system by choosing the law $\P^u$, equivalent to a reference probability $\P$, under which $B^u:=(B^u_s)_{s\le T}$ is a Brownian motion. On the other hand, at her/his convenience, the controller chooses also the time $\tau$ to stop controlling the system. As a result this incurs a payoff, which is a reward, $J(u,\tau)$ given by: 
\begin{equation}\label{pay1}
J(u,\tau):=\E^u\left[\int_{0}^{\tau}\Gamma(s,x,u_s)ds+L_\tau\ind_{\{\tau<T\}}+g(x)\ind_{\{\tau=T\}}\right].
\end{equation} Here: i)  $\Gamma$ is the instantaneous reward; ii) $L_\t$ is the reward if the controller decides  to stop at $\tau$ before the terminal time $T$; iii) $g$ is the reward at the terminal $T$. 

This problem is termed of control-stopping (or mixed) type because it combines control and stopping. It has been considered in several papers including \cite{BY, nicoleaspects,hlp00, KZ}. There are at least two methods to tackle this problem. One is based on the martingale approach (\cite{nicoleaspects, KZ}) and the other one uses the notion of backward stochastic differential equations (\cite{BY,hlp00}). However in all those works there are technical restrictions on the data ($f$, $\G$, $L$, $g$, etc.) which define this problem. Actually, in \cite{nicoleaspects, KZ}, $f$ is supposed of linear growth w.r.t. $x$ and $\G$, $L$, $g$ bounded while in \cite{BY} this latter assumption of boundedness is relaxed to linear growth. On the other hand, in \cite{hlp00}, $f$ is supposed bounded while the other functions can have an arbitrary polynomial growth w.r.t. $x$. Therefore the main objective of this work is to unify those frameworks, i.e., to consider the control-stopping problem when $f$, on one hand, and $\G$, $L$, $g$, on the other hand, are of linear and polynomial growths respectively. The approach is based on reflected BSDEs. 

To deal with the control-stopping in our general framework we are led to study the following reflected BSDE \eqref{EK1intro} and to prove that it has a solution:  
 \begin{equation}\label{EK1intro}
 \left\{\begin{array}{l}
 Y_t = g(x) + \int_t^T
 H^*(s,x,Z_s) ds +K_T-K_t- \int_t^T Z_s dB_s,\,\,t\le T;\,
\\\\
 L\leq Y ~\mx{and}~\int_{0}^{T}\left(Y_s-L_s\right)dK_s=0,
 \end{array}\right.
 \end{equation}
where $$\bal H^*(s,x,z)=\sup_{u\in\cal{A}}H(s,x,z,u)\mx{ with}
\\\\\qquad H(t,x,z,u)=z\sigma^{-1}(t,x)f(t,x,u)+\Gamma(t,x,u).\ea$$ The function $H$ is the Hamiltonian of the problem and $\ca$ the set of values of the controls. 

First let us notice that if $\sigma^{-1}f$ is bounded then obviously the 
function $H^*$ is Lipschitz w.r.t. $z$ and consequently the existence-uniqueness of a solution for the above reflected BSDE \eqref{EK1intro} is obtained from the classical result by El-Karoui \cite{EKPPQ}. Now in the case when $\sigma^{-1}f$ is not bounded and satisfies a linear growth condition only, $H^*$ verifes the following property which is called the stochastic Lipschitz condition: $\fr$ $z,z'\in \R^d$,
\begin{equation}\lb{croipolyhstar}
|H^*(t,x,z)-H^*(t,x,z)|\le C(1+\|x\|_t)|z-z'|\,\,(\|x\|_t=\sup_{s\le t}|x_s|).
\end{equation}
which means that the reflected BSDE \eqref{EK1intro} is not standard and of stochastic Lipschitz type. These latter RBSDEs are already considered in some papers including \cite{moth, W} (for non reflected ones, see \cite{BK}). However the results of these works do not allow to deduce satisfactorily the existence of a solution to \eqref{EK1intro} since they have been stated in frameworks which do not fit completely to ours. Indeed, had we applied those results, we would have been led to assume restrictive conditions on the horizon $T$ of the problem which, among other conditions, should be small. Therefore the first task is to show that the reflected BSDE \eqref{EK1intro} has a solution for arbitrary finite horizon $T$. Under condition \eqref{croipolyhstar} on $H^*$ and the polynomial growth of $\G$, $g$ and $L$ with respect to $x$, we show that the reflected BSDE \eqref{EK1intro} has a unique solution. In this proof, one point plays a crucial role which is the existence of a constant $p>1$ such that
$$
\E[(\frac{d\P^u}{d\P})^p]<\infty.$$
Further the link with the control-stopping problem is stated, an optimal pair $(u^*,\tau^*)$ is exhibited and finally the associated  Hamilton-Jacobi-Bellman equation, in the Markovian framework, is studied. 

The paper is organized as follows. In Section 2, we formulate the problem and consider a specific reflected BSDE with stochastic Lipschitz coefficient whose solution provides appropriate estimates for the solution of the reflected BSDE \eqref{EK1intro}. In Section 3, we show that equation \eqref{EK1intro} has a unique solution $(Y,Z,K)$. This solution is constructed as a limit (twice) of an approximating scheme obtained by truncating $H^*$  twice. Namely we show that 
$$
Y=\lim_{m\rw \infty}\lim_{n\rw \infty} Y^{n,m}$$
where $(Y^{n,m},Z^{n,m},K^{n,m})$ is the solution of the standard reflected BSDE associated with $(H^{*n,m},g(x),L)$ with 
$$
H^{*n,m}(t,x,z):=\max\{H^*(t,x,z),0\}\ind_{\{\|x\|_t\le n\}}
-\max\{-H^*(t,x,z),0\}\ind_{\{\|x\|_t\le m\}}$$which is Lipschitz w.r.t. $z$ since $H^*$ verifies \eqref{croipolyhstar}. Another property which plays an important role in this construction is the comparison of solutions of reflected BSDEs. We then show that $Y_0$ is nothing else but the optimal payoff of the control-stopping problem. The optimal pair of control and stopping time is exhibited. 

In Section 4, we consider the Markovian framework of the mixed control problem, i.e., roughly speaking, for $\txt$,  the dynamics of the controlled system is given by:
\begin{equation}
dx_s^{t,x}=f(t,x^{t,x}_s,u_s)ds+\sigma(s,x^{t,x}_s)dB^u_s, \,\,t\in[0,T] \mx{ and } x^{t,x}_t=x \in \R^d\mx{ fixed}.
\end{equation}
The payoff $J_t(u,\tau)$ on the time interval $[t,T]$,  is defined by: \begin{equation}\label{pay12}
J_t(u,\tau):=\E^u\left[\int_{t}^{\tau}\Gamma(s,x^{t,x}_s,u_s)ds+h(\tau,x^{t,x}_\tau)\ind_{\{\tau<T\}}+g(x^{t,x}_T)\ind_{\{\tau=T\}}\right].
\end{equation} We show that the deterministic function 
$$
u(t,x):=Y^{t,x}_t=\sup_{(u,\t)}J_t(u,\tau),$$ 
where $(Y^{t,x},Z^{t,x},K^{t,x})$ is the unique solution of the reflected BSDE \eqref{EK1intro} in this Markovian framework, is the unique viscosity solution of the Hamilton-Jacobi-Bellman associated with the control-stopping problem, i.e., 
\begin{equation}\label{hjbintro}
\left\{
\begin{array}{l}
\min\left[u(t,x)-h(t,x),
-\partial_t u(t,x)-\cl u(t,x)-H^*(t,x,\nabla_xu(t,x)\sigma(t,x))\right]=0, \,\,(t,x)\in[0,T[ \times \R ^d;\\\\
u(T,x)=g(x),\, x\in\R^d.
\end{array}
\right.
\end{equation}
Moreover $u$ is continuous and of polynomial growth. The main difficulty is to show continuity of $\bar u^m(t,x)=\lim_{n\rw \infty}Y^{t,x,n,m}_t$. Due to the polynomial growths of $\G$, $h$ and $g$, this continuity cannot be obtained by the usual characterization by means of Snell envelope of processes. To overcome this difficulty we have shown that the comparison principle holds for the partial differential equation (PDE for short) associated with the RBSDE verified by $Y^{t,x,m}=\lim_{n\rw \infty}Y^{t,x,n,m}$ (see \eqref{vis-m}) and that $\bar u^m(t,x)$ is a solution. Consequently it is continuous and unique in the classe of functions with polynomial growth.  Finally, similarly, we deduce that $u(t,x)$ is continuous and is the unique solution of \eqref{hjbintro} in the classe of functions with polynomial growth.\qed
\section{Formulation of the problem. Study of a specific reflected BSDE}

Let $\Omega=\cc\left([0,T];\R^d\right)$ be the space of $\R^d$-valued continuous function on $[0,T]$ endowed with the metric of uniform convergence on $[0,T]$.
Denote by $\F$ the Borel $\sigma$-field over $\Omega$. Next for $\o\in \O$ and $t\le T$, let us set $||\omega||_t:=\sup_{0\leq s\leq t}|\omega_s|$. Let $x:=(x_s)_{s\le T}$ be the coordinate process on $\O$, i.e., $x_s(\o)=\o(s)$ and denote by $(\F^0_t:=\sigma(x_s,s\leq t))_{t\in [0,T]}$, the filtration on $\O$ generated by $x$.\\

Let $\sigma$ be a function from $[0,T]\times\Omega$ into $\R^{d \times d}$ such that:\\
\begin{enumerate}
\item[A1)]$\sigma$ is $\F^0_t$-progressively measurable.
\item[A2)] There exists a constant $C>0$ such that:
\begin{enumerate}
\item[a)] For every $t\in [0,T]$ and $\omega,\omega'\in \Omega,|\sigma(t,\omega)-\sigma(t,\omega')|\leq C||\omega-\omega'||_t$,
\item[b)]$\sigma$ is bounded, invertible and its inverse $\sigma^{-1}$ is bounded.
\end{enumerate}
\end{enumerate}
Let $\P$ be a probability measure on $\Omega$ such that $(\Omega, \P)$ carries a $d$-dimensional Brownian motion $(B_t)_{0\leq t \leq T}$ and for any  $x_0 \in \R^d$, the process $(x_t)_{0\leq t \leq T}$ is the unique solution of the
following stochastic differential equation:
\begin{equation}\ba{l}
x_t=x_0+ \int_0^t\sigma(s,x)dB_s,\quad t\in [0,T].\ea
\end{equation}
Such a pair $(\P, B)$  exists thanks to Proposition 4.6 in (\cite{KS}, pp.315) since $\sigma$ satisfies (A2). Moreover, for any $p \geq 1$,
\begin{equation}\label{estimx}
\E\left[\|x\|^p_T\right]\leq C_p
\end{equation}
where $C_p$ depends only on $p$, $T$, the initial value $x_0$ and the linear growth constant of $\sigma$ (see
\cite{KS}, pp. 306). Again, since $\sigma$ satisfies (A2), $\F^0_t$ is the same as $\sigma \{B_s,s \leq t\}$ for any $t \leq T$.
We denote by $\F:= (\F_t)_{0\leq t\leq T}$ the completion of $(\F^0_t)_{t\leq T}$ with the $\P$-null sets of $\Omega$. Finally let $\cp$ be the $\F$-progressively measurable $\sigma$-algebra on $\O\times [0,T]$.\\

Let us introduce some notations. In the following, for $q,k\geq 1$, we denote by:
\bs

\ni i) $ L^{q}:=$\{$\xi$ is an $\F_T$-measurable random variable s.t $\E[|\xi|^{q}] <+\infty$\}.\\
\ni ii) $\H^{q,k}:=\{Z:=(Z_s)_{s\le T}$, $\cp$-measurable, $\R^k$-valued process s.t. $\E\left[\left(\int_{0}^{T}|Z_s|^2\right)^{\frac{q}{2}}\right]<+\infty\}$.\\
\ni iii) $\cs^q$:=\{$Y:=(Y_t)_{t\leq T}$, $\cp$-measurable, $\R$-valued and continuous process s.t. $\E\left[ \sup_{0\leq t\leq T}|Y_t|^q\right]<+\infty$ \}.\\
\ni iv) By convention $\inf\{\emptyset\}=+\infty$.\\
\ms

Let us now consider the following functions $g$, $\varphi$ and obstacle $(L_t)_{t\leq T}$.

\begin{enumerate}
\item[(i)] \ \ $g: \cc \rightarrow \R$ is a Borel measurable function of polynomial growth, i.e., there exist constants $C$ and $p$ such that:
\begin{center}
$ |g(x)|\leq C(1+{\|x\|_T^p}),\quad \forall x\in\cc.$
\end{center}
\item[(ii)] The process $(L_t)_{t\le T}$ is $\cp$-measurable, continuous, $\R$-valued. Moreover there exist  positive constants $C$ and $p$ such that $\forall t\in[0,T]$,
\begin{equation}\lb{grcdtl}|L_t|\leq C(1+{\|x\|_t^p}).\end{equation}
We moreover assume that $\P$-a.s. $L_T\le g(x)$.

\item[(iii)] The function $\varphi$ is defined as follows: For some $p\ge 1$, for any $(t, x, z)\in[0,T]\times \cc \times \R^d$,
$$ \varphi(t,x,z)=c(1+\|x\|_t)|z|+c(1+{\|x\|_t^p}).$$
\end{enumerate}

We now introduce a property of exponential martingales which plays an important role in this work. Let $\vartheta$ be a function from
$[0,T]\times \cc$ into $\R^d$ which is $\cp$-measurable. For $t\le T$ we set:
$$\zeta_t=e^{\int_0^t\vartheta(s,\omega)\sigma^{-1}(s,\omega)dB_s-\frac{1}{2}\int_0^t|
\vartheta(s,\omega)\sigma^{-1}(s,\omega)|^2ds}.$$
We then have:
\begin{lemma}\label{lemmeproba}${}$\\i) Assume that $\vartheta$ is bounded. Then
$(\zeta_t)_{t\le T}$ is a martingale such that for any $\ell \ge 1$,
$$
\E[(\zeta_T)^\ell]\le C_{\ell,\vartheta, T}
$$
where $C_{\ell,\vartheta,T}$ is a constant which depends only on $T$, $\ell$ and the constant of boundedness of $\vartheta \sigma^{-1}$.
\ms

\ni ii) Assume that $\vartheta$ is of linear growth, i.e.,
for any $t\le T$, $|\vartheta(t,\omega)|\le c_{\vartheta}(1+|\omega|_t)$. Then there exist constants $q>1$ and $\varpi$ which depend only on $T$, $\sigma$, $\sigma^{-1}$ and $c_{\vartheta}$ such that
\begin{equation}\label{probhaus}
\E[(\zeta_T)^q]\le \varpi.
\end{equation}
Moreover $(\zeta_t)_{t\le T}$ is a martingale.
\ms

\ni iii) In both cases, the measure $Q$ such that $dQ=\zeta_T.d\P$ is a probability on $(\Omega,\F)$.
\end{lemma}
\ni \textbf{Proof}: Point i) is classical and based on the fact that when $\vartheta$ is bounded, the familly of $\zeta$'s are martingales. The proof of point ii) is given in (\cite{UG}, Equation (2.8), pp. 3). Finally point iii) is obvious since $\zeta$ is a positive martingale. \qed
\bigskip

Next let $\Phi$ be a function from $[0,T]\times \cc\ \times\R^d$ into $\R$ which is $\cp\otimes {\bf B}(\R^d)$-measurable. First we define the notion of a solution of the reflected BSDE associated
with a triplet $(g, \Phi, L)$ which we consider throughout this paper.
\begin{definition}\label{defiLpsol}A triplet $(Y,Z,K):=(Y_t,Z_t,K_t)_{t\le T}$ of $\cp$-measurable and $\R^{1+d+1}$-valued processes  is a
solution of the reflected BSDE associated with the lower
reflecting barrier $L$, the terminal  value $g(x)$ and generator $\Phi$
if the followings hold:
\ms

\ni i) $Y$ is continuous, $K$ is continuous non-decreasing ($K_0=0$), $\P$-a.s. $(Z_t(\omega))_{t\leq T}$ is $dt$-square integrable;\\
ii) $Y_t = g(x) + \int_t^T \Phi(s,x,Z_s) ds + K_T - K_t - \int_t^T Z_s dB_s$, $0 \leq t \leq T$;\\
iii) $Y_t \geq L_t$, $0 \leq t \leq T$ and $\int_0^T (Y_s - L_s) dK_s =0$.
\end{definition}

We are going to show that the reflected BSDE associated with $(\varphi, L,g)$ has a solution which also satisfies other integrability properties. As it is mentionned previously, $\varphi$ is not a standard generator which does not enter neither in the framework of \cite{EKPPQ} nor in the one of \cite{BK}.
\begin{proposition}\label{unique}
There exist $\cp$-measurable processes $(Y, Z, K)$ valued in $\R^{1+d+1}$ and a stationary non-decreasing sequence of
stopping times $(\tau_k)_{k\geq 1}$ such that:
\medskip

\ni i) $(Y,Z,K)$ is a solution for the BSDE associated with $(\varphi, L, g)$.

\ni ii) For any constant $\gm\ge 1$ and $\t$ a stopping time valued in $[0,T]$,
\begin{equation}\lb{ygamma}
\E[|Y_\t|^\gm]<+\infty.
\end{equation}
\ni iii) For any $k\geq 1$,
\begin{equation}\label{estimloc}\ba{l}
\E[\sup_{s\leq T}|Y_{s\wedge \t_k}|^\gm+K^\gm_{\t_k}+\int_0^{\t_k}|Z_s|^2ds]<+\infty,\ea\end{equation}
where $\t_k$ depends only on $g$, $L$ and $x$.
\end{proposition}
\ni \textbf{Proof}: For $n\geq 0$ let us set
$$\varphi^n(t,x,z)=[c(1+\|x\|_t)\wedge n]|z|+c(1+{\|x\|_t^p}).$$
Then $\varphi_n$ is Lipschitz w.r.t $z$. Therefore by El-Karoui et al.'s result in \cite{EKPPQ}, there exists a triplet of processes
$(Y^n,Z^n,K^n)\in \cs^2\times \H^{2,d}\times \cs^2 $ ($K^n$ is non-decreasing and $K^n_0=0$) such that :
\begin{equation}\label{impn}
\left\{\begin{array}{l}
Y^n_t = g(x) + \int_t^T
 \varphi^n(s,x,Z^n_s)ds +K^n_T-K^n_t- \int_t^T Z^n_s dB_s,\quad
t\in[0,T]\,  ;\\
L_t \leq Y_t^n,\, \fr t\in[0,T] \mbox{ and }\int_{0}^{T}\left(Y_s^n-L_s\right)dK^n_s=0.
\end{array}\right.
\end{equation}
Next by the comparison Theorem (see e.g. \cite{EKPPQ}) we have, for any $n\ge 0$, $Y^n \leq Y^{n+1}$ and $dK^{n+1}\leq dK^n$ since $\varphi^n\leq \varphi^{n+1}$ and the barrier $L$ is fixed.

We now show the following Lemmas (which are steps forward in the proof of Proposition \ref{unique}) related to  estimations of the processes $Y^n, n\ge 0$, and the convergence of the sequence $(Y_n)_{n\ge 0}$.

\begin{lemma} \label{lemme1}There exists a $\cp$-measurable RCLL (for right continuous with left limits) process $Y:=(Y_t)_{t\leq T}$, $\R$-valued such that $\P$-a.s. for any $t\leq T$, $Y^n_t\nearrow Y_t$. Moreover for any $\gm\ge 1$ and any stopping times $\t\in[0,T]$, $$\E\left[|Y_\t|^\gm \right]\leq c
$$where $c$ is a constant independent of $\t$.
\end{lemma}
\ni \textbf{Proof}: Let $\P^n$ be the probability, equivalent to $\P$, defined as follows: $$d\P^n=L^n_T d\P$$where for any $t\le T$,
$$\ba{l}
L^n_t:=\exp\{\int_0^t \{c(1+||x||_s)\wedge n\}\ell(Z^n_s) dB_s-
\frac{1}{2}\int_0^t\|\{c(1+\|x\|_s)\wedge n\}\ell(Z^n_s)\|^2 ds\}\ea
$$
with $\ell$ is a bounded measurable function such that $\ell(z).z=|z|$, $\fr z=(z_i)_{i=1,\dots,d}\in \R^d$. Such a function $\ell$ exists and it is enough to take $\ell(z)=(\ell_i(z))_{i=1,\dots,d}$ where for $i=1,\dots,d-1$\\
$\ell_i(z)=(\sqrt{z_i^2+...+z_d^2}-\sqrt{z_{i+1}^2+...+z_d^2})z_i^{-1}\ind_{\{z_i\neq 0\}}$ and 
$\ell_d(z)=|z_d|z_d^{-1}\ind_{\{z_d\neq 0\}}$. Then by Girsanov's Theorem, the process $(B^n_t:=B_t-\int_0^t
 \{c(1+\|x\|_s)\wedge n\}\ell(Z^n_s)ds)_{t\ge 0}$ is a Brownian motion under $\P^n$ and the triplet $(Y^n,Z^n,K^n)$ verifies: For any $t\in[0,T]$,
$$\left\{\begin{array}{l}
Y^n_t = g(x) + \int_t^T
 c(1+{\|x\|_s^p})ds +K^n_T-K^n_t- \int_t^T Z^n_s dB^n_s\,\,;\\
L_t \leq Y_t^n  \mbox{ and }\int_{0}^{T}\left(Y_s^n-L_s\right)dK^n_s=0.
\end{array}\right.$$
Therefore for any stopping time $\t\leq T$, $\P$-a.s., we have:
\begin{eqnarray*}
Y^n_\t = \mbox{esssup}_{\sigma \geq \t}\E^{\P^n}\left[\int_\t^\sigma c(1+\|x\|^p_s)ds+L_\sigma \ind_{\{\sigma < T\}}+g(x_T) \ind_{\{\sigma=T\}}|\F_\t\right]
\end{eqnarray*}
(one can see \cite{EKPPQ} for this characterization). Next by using polynomial growth of $g$ and $L$ and the fact that $\P$ and $\P^n$ are equivalent, we deduce that:
\begin{eqnarray*}
|Y^n_\t| && \le \E^{\P^n}\left[\int_\t^Tc(1+\|x\|^p_s)ds+2c(1+\|x\|^p_T)|\F_\t\right]\\
&&\leq c(1+\E^{\P^n}\left[\|x\|^p_T|\F_\t\right]).
\end{eqnarray*}
Next let $\gm\geq 1$. Then by conditional Jensen's inequality we have:
\begin{align}\label{estimyn}
|Y^n_\t|^\gm &\leq c(1+\E^{\P^n}\left[\|x\|^{\gm p}_T|\F_\t\right])\\ \nonumber
&=c(1+\E\left[L^n_{\t,T}\|x\|^{\gm p}_T|\F_\t\right])
\end{align}
where for any $t\leq T$,
\begin{align}&L^n_{t,T}:={L^n_T}\div {L^n_t}\nn\\&\,\,=\exp\{\int_t^T \{c(1+||x||_s)\wedge n\}\ell(Z^n_s) dB_s-
\frac{1}{2}\int_t^T\|\{c(1+\|x\|_s)\wedge n\}\ell(Z^n_s)\|^2 ds\}\nn\\&
=\exp\{\int_0^T1_{\{t\le s\le T\}}\left( \{c(1+||x||_s)\wedge n\}\ell(Z^n_s) dB_s-
\frac{1}{2}\|\{c(1+\|x\|_s)\wedge n\}\ell(Z^n_s)\|^2 ds \right)\}
.\end{align} Now by Lemma \ref{lemmeproba}, there exist constants $c$ and $q > 1$, which do not depend on $n$ and $\t$ such that $\E[( L^n_{\t,T})^q]\leq c$. Next by Young's inequality we get from \eqref{estimyn}:
\begin{eqnarray*}
|Y^n_\t|^\gm\le
c(1+\E\left[q^{-1}(L^n_{\t,T})^q+\check{q}^{-1}\|x\|^{\gm p\check{q}}_T|\F_\t\right]),
\end{eqnarray*}
where $q^{-1}+\check q^{-1}=1$. Taking now into account of \eqref{estimx} we deduce that:
\begin{equation}\label{eqyng}\fr n\ge 0,\,\,\E[|Y^n_\t|^\gm]\leq c\end{equation}
where $c$ is a constant which does not depend on $n$. As for any $n\ge 0$, $Y^n\le Y^{n+1}$ then $\P$-a.s.
for any $t\leq T$, $Y^n_t\rightarrow_n Y_t=\liminf_{n\rightarrow +\infty}Y^n_t$. Therefore $\P$-a.s. $Y^n_\t\rightarrow_n Y_\t$ and by Fatou's Lemma and \eqref{eqyng} we have $\E[|Y_\t|^\gm]\leq c$.

It remains to show that $Y$ is RCLL. But this a direct consequence of the fact that $Y^n$ is a continuous supermartingale which converges increasingly and pointwise to $Y$ and then the limiting process $Y$ is also RCLL (see \cite{DMM}, pp.86).
\ms

Next let $\tilde \theta=|x_0|+|L_0|+|Y^0_0|+|Y_0|$ ($\tilde \theta$ is a constant) and for $k\ge 1$ let us define the sequence of stopping times $(\t_k)_{k\ge 1}$ by:
$$\tau_k:=\inf\{t\geq 0, |Y_t|+\|x\|_t+|L_t|+|Y^0_t|\ge \tilde \theta +k \}\wedge T.$$
The sequence of stopping times $(\tau_k)_{k\geq 1}$
is non-decreasing, of stationary type (i.e. constant after some rank $k_0(\o)$) converging to T
since the process $Y$ is RCLL and $Y^0$, $x$ and $L$ are continuous. Moreover for any $k\ge 1$,
$$ max\{\sup_{t\leq \tau_k}|L_t|,\sup_{t\leq \tau_k}|Y_t|,\sup_{t\leq \tau_k}|Y^n_t|\} \leq \tilde \theta + k:=\tilde \theta_k.$$
Next we have the following result:
\begin{lemma} \label{lemme2}${}$\\
\ni i) The process $Y$ is continuous.

\ni ii) There exist $\cp$-measurable processes $(K, Z)$ valued in $\R^{1+d}$ such that $(Y,Z,K)$ is a solution of the reflected BSDE associated with $(\varphi, g,L)$ and verifies \eqref{estimloc}.\end{lemma}
\ni \textbf{Proof}: For any $k\geq 1$ and $n\geq 0$ we have: $\forall t\in[0,T]$
$$\left\{ \begin{array}{l}Y^n_{t\wedge\tau_k} = Y^n_{\tau_k} + \int_{t\wedge\tau_k}^{\tau_k}
 \varphi^n(s,x,Z^n_s)ds +K^n_{\tau_k}-K^n_{t\wedge\tau_k}- \int_{t\wedge\tau_k}^{\tau_k} Z^n_s dB_s;
\\\\
L_{t\wedge\tau_k}\leq Y^n_{t\wedge\tau_k}\mbox{ and }\int_{0}^{\t_k}\left(Y^n_s-L_s\right)dK^n_s=0.
\end{array}\right.$$
By using the It\^o formula with $(Y^n_{t\wedge\tau_k})^2$ and taking into account of \eqref{estimx}, we classicaly deduce the existence of a constant $C_k$, which depends on $k$, such that uniformly on $n$, we have
\begin{equation}\label{estimznk}
\E\left[\int_0^{\t_k}|Z^n_s|^2ds\right]\le C_k
\end{equation}
since $|Y^n_{t\wedge\tau_k}|\le \tilde \theta_k$, for any $t\le T$, and $Y^0\le Y^n\le Y$.

Next once more by It\^{o}'s formula  we obtain: $\forall t\in[0,T]$,
\begin{align}
(Y^n_{t\wedge\tau_k}-Y^m_{t\wedge\tau_k})^{2}&=(Y^n_{\tau_k}-Y^m_{\tau_k})^{2}+2\int_{t\wedge\tau_k}^{\tau_k}(Y^n_s-Y^m_s)(\varphi^n(s,x,Z_s^n)-\varphi^m(s,x,Z_s^m))ds\nn\\ \nonumber
& \qquad +2\int_{t\wedge\tau_k}^{\tau_k}(Y^n_s-Y^m_s)d(K^n_s-K^m_s)-\int_{t\wedge\tau_k}^{\tau_k}|Z^n_s-Z^m_s|^2
ds\\&\nn \qquad -2\int_{t\wedge\tau_k}^{\tau_k}(Y^n_s-Y^m_s)(Z^n_s-Z^m_s)dB_s.
\end{align}
The definition of $\varphi^n$ and the fact that $(Y^n_s-Y^m_s)d(K^n_s-K^m_s)\le 0$ imply that: $\forall t\le T$,
\begin{align}\label{deducsup}
&(Y^n_{t\wedge\tau_k}-Y^m_{t\wedge\tau_k})^{2}+\int_{t\wedge\tau_k}^{\tau_k}|Z^n_s-Z^m_s|^2ds\le (Y^n_{\tau_k}-Y^m_{\tau_k})^{2}\\&\nn\qquad\qquad\qquad +2\int_{t\wedge\tau_k}^{\tau_k}|Y^n_s-Y^m_s|C(1+\|x\|_s)(|Z_s^n|+|Z_s^m|)ds-2\int_{t\wedge\tau_k}^{\tau_k}(Y^n_s-Y^m_s)(Z^n_s-Z^m_s)dB_s.\end{align} Take expectation on both hand-sides to deduce that
\begin{align}
&\nn \E[\int_{t\wedge\tau_k}^{\tau_k}|Z^n_s-Z^m_s|^2ds]\le \E[(Y^n_{\tau_k}-Y^m_{\tau_k})^{2}]+2\E[\int_{t\wedge\tau_k}^{\tau_k}|Y^n_s-Y^m_s|C(1+\|x\|_s)(|Z_s^n|+|Z_s^m|)ds].\end{align}
Now the definition of $\t_k$, estimate \eqref{estimznk} and the Cauchy-Schwarz inequality yield:
\begin{equation}\lb{convznmk}\E[\int_{t\wedge\tau_k}^{\tau_k}|Z^n_s-Z^m_s|^2ds]\rightarrow 0 \mbox{ as }n,m\rightarrow +\infty.\end{equation}
Consequently the sequence $((Z_s^n1_{\{0\le s\le \t_k\}})_{s\leq T})_{n\ge 0}$ is of Cauchy type in $\H^{2,d}$ and then there exists a process $Z_k$ which belongs to $\H^{2,d}$ such that  $((Z_s^n1_{\{0\le s\le \t_k\}})_{s\leq T})_{n\ge 0}\rightarrow_n Z_k(s)$ in $\H^{2,d}$.

Next going back to \eqref{deducsup}, take the supremum over $t$, make use of Burkholder-Davis-Gundy (see e.g. \cite{KS, [RY]} and BDG for short) inequality and finally take the expectation to deduce that
\begin{equation}\lb{convynm}
\E[\sup_{s\leq T}|Y^n_{s\wedge \t_k}-Y^m_{s\wedge \t_k}|^2]\rightarrow 0 \mbox{ as }n,m\rightarrow +\infty.
\end{equation}
As $(Y^n)_{n\ge 0}\rightarrow Y$ then the process $(Y_{s\wedge \t_k})_{s\ge 0}$ is continuous for any $k\ge 1$ and
\begin{equation}\lb{convnk}
\E[\sup_{s\leq T}|Y^n_{s\wedge \t_k}-Y_{s\wedge \t_k}|^2]\rightarrow 0 \mbox{ as }n\rightarrow +\infty.
\end{equation}
But the sequence $(\t_k)_{k\ge 1}$ is of stationary type, therefore the process $Y$ is continuous.

Now by \eqref{estimznk} and \eqref{convznmk}, for any $k\ge 1$, the sequence of processes
$((\varphi^n(s,x_s,Z^n_s)1_{\{s\le \t_k\}})_{s\le T})_{n\ge 0}$ converges in $\H^{2,d}$ to
$(\varphi(s,x,Z_k(s))1_{\{s\le \t_k\}})_{s\le T})$. Therefore if we set, for $k\geq 1$ and $t\le T$,
$$
K_k(t)=Y_0-Y_{t\wedge \tau_k}- \int_0^{t\wedge\tau_k}\varphi(s,x,Z_k(s))ds +\int^{t\wedge\tau_k}_0 Z_k(s) dB_s
$$
we obtain that
\begin{equation}\lb{convkn}
\E[\sup_{s\leq T}|K^n_{s\wedge \t_k}-K_k(s)|^2]\rightarrow 0 \mbox{ as }n\rightarrow +\infty.
\end{equation}
Finally the uniform convergences \eqref{convkn} and \eqref{convnk} imply that, 
in view of Helly's Convergence Theorem (see \cite{KF70}, pp. 370), 
$$
\int_0^{\t_k}(Y_s-L_s)dK_k(s)=0.
$$
It means that for any $k\ge 1$ we have: $\fr t\leq T$,
\begin{equation}\lb{eqk1}\left\{ \begin{array}{l}Y_{t\wedge\tau_k} = Y_{\tau_k} + \int_{t\wedge\tau_k}^{\tau_k}
 \varphi(s,x,Z_k(s))ds +K_k({\tau_k})-K_k({t\wedge\tau_k})- \int_{t\wedge\tau_k}^{\tau_k} Z_k(s)dB_s;
\\\\
L_{t\wedge\tau_k}\leq Y_{t\wedge\tau_k}\mbox{ and }\int_{0}^{\t_k}\left(Y_s-L_s\right)dK_k(s)=0.
\end{array}\right.\end{equation}
Take now the reflected BSDE satisfied by $(Y,Z_{k+1},K_{k+1})$ on $[0,\t_k]$ (since $\t_k\le \t_{k+1}$) yields: for any $t\leq T$,
\begin{equation}\lb{eqk2}\left\{ \begin{array}{l}Y_{t\wedge\tau_k} = Y_{\tau_k} + \int_{t\wedge\tau_k}^{\tau_{k}}
 \varphi(s,x,Z_{k+1}(s))ds +K_{k+1}({\tau_k})-K_{k+1}({t\wedge\tau_k})- \int_{t\wedge\tau_k}^{\tau_k} Z_{k+1}(s)dB_s;
\\\\
L_{t\wedge\tau_k}\leq Y_{t\wedge\tau_k}\mbox{ and }\int_{0}^{\t_k}\left(Y_s-L_s\right)dK_{k+1}(s)=0.
\end{array}\right.\end{equation}
By uniqueness (which holds since on $[0,\t_{k}]$, $\varphi$ is Lipschitz w.r.t $z$), we have for any $k\ge 1$ :
$$
Z_k(s)1_{\{s\le \t_k\}}=Z_{k+1}(s)1_{\{s\le \t_k\}} \,\,ds\otimes d\P-a.e \mbox{ and }K_{k+1}({s\wedge\tau_k})=K_{k}({s\wedge\tau_k}),\,\,\fr s\le T.
$$
Then let us define (by concatenation) the processes $Z$ and $K$ as follows: $\forall t\le T$,
$$
Z_t=Z_1(t)1_{\{t\le \t_1\}}+\sum_{k\ge 2}Z_k(t)1_{\{\t_{k-1}< t\le \t_{k}\}}$$ and $$K_t=\left\{\begin{array}{l}
K_1(t) \mbox{ if }t\leq \t_1;\\\\
(K_{k+1}(t)-K_{k+1}(\t_k)) +K(\t_k)\mbox{ for }\t_k<t\leq \t_{k+1}, \,\,k\geq 1.
\end{array}\right.
$$
Note that the processes $Z$ and $K$ are well-defined since the sequence ($\t_k)_{k\ge 1}$ is of stationary type. On the other hand, $K$ is continuous non-decreasing and $\P$-a.s., $K_T(\o)<+\infty$ and $(Z_s(\o))_{s\le T}$ is $ds$-square integrable. Finally for any $k\ge 1$, it holds:
\begin{equation}\lb{eqk1deux}\left\{ \begin{array}{l}Y_{t\wedge\tau_k} = Y_{\tau_k} + \int_{t\wedge\tau_k}^{\tau_k}
 \varphi(s,x,Z_s)ds +K_{\tau_k}-K_{t\wedge\tau_k}- \int_{t\wedge\tau_k}^{\tau_k} Z_s dB_s;
\\\\
L_{t\wedge\tau_k}\leq Y_{t\wedge\tau_k}\mbox{ and }\int_{0}^{\t_k}\left(Y_s-L_s\right)dK_s=0.
\end{array}\right.\end{equation}
Take now $k$ great enough and since once more ($\t_k)_{k\ge 1}$ is of stationary type to obtain: $\fr t\le T$,
\begin{equation}\lb{eqk1trois}\left\{ \begin{array}{l}Y_{t} = g(x)+ \int_{t}^{T}
 \varphi(s,x,Z(s))ds +K_T-K_t- \int_{t}^{T} Z_s dB_s;
\\\\
L_{t}\leq Y_{t}\mbox{ and }\int_{0}^{T}\left(Y_s-L_s\right)dK_s=0
\end{array}\right.\end{equation} which completes the proof. \qed
\section{Study of the stochastic mixed control problem}
Let $\ca$ be a compact metric space and let ${\cal U}$ be the space of
${\cal P}$-measurable processes $u:=(u_t)_{t\leq T}$ with values in $\ca$. We first introduce the objects and assumptions which we need in this section.
\begin{enumerate}
\item[(1)] $f:[0,T]\times \Omega \times \ca\rightarrow \R^d$ is a function such that:
\begin{enumerate}
\item[i)] For each $a\in \ca$, the function $(t,\o)\rightarrow f(t,\o,a)$ is $\cp$-measurable.
\item[ii)] For each $(t,\o)$, the mapping $a\rightarrow f(t,\o,a)$ is  continuous.
\item[iii)] $f$ is of linear growth, i.e., there exists a real constant $C>0$ such that:
\begin{equation}\label{generateur-cond1}
|f(t,\o,a)|\leq C(1+||\o||_t),\, \forall\, 0\leq t \leq T,\,\o\in
\Omega,\,a\in \ca.
\end{equation}
\end{enumerate}
\item[(2)] $\Gamma:[0,T]\times \Omega \times \ca\rightarrow \R$ is a $\cp \otimes {\bf B}(\ca)$-measurable function such that for each $(t,\o)$ the mapping $a\rightarrow \Gamma(t,\o,a)$
is  continuous. In addition there exist positive constants $C>0$ and $p$ such that:
\begin{equation}\label{hgenerateur1}
|\Gamma(t,\o,a)|\leq C(1+\|\o\|_t^p),\, \forall \,t \leq T,\,w\in
\Omega,\,a\in \ca.
\end{equation}
\end{enumerate}

Next let $u:={\left(u_t\right)}_{t\leq T}$ be an adapted $\ca$-valued stochastic process and $\P^u$ the probability on $\cc$ (which is actually one) equivalent to $\P$ such that:
$$d\P^u= M_Td\P$$where for any $t\le T$,
$$
M_t:=e^{\int_0^t\sigma^{-1}(s,x)f(s,x,u_s)dB_s-
\frac{1}{2}\int_0^t\|\sigma^{-1}(s,x)f(s,x,u_s)\|^2ds}.
$$
Under $\P^u$ the process ${\left(x_t\right)}_{t\leq T}$ verifies:
\begin{equation}\label{eds-contr}\begin{array}{l}
x_t=x_0+\int_0^tf(s,x,u_s)ds+\int_0^t\sigma(s,x)dB^u_s,\quad  t
\leq T,\end{array}
\end{equation}
where
\begin{equation}\label{chang-prob}\begin{array}{l}
B^u_t:=B_t-\int_0^t\sigma^{-1}(s,x)f(s,x,u_s)ds,\quad  t \leq T,\end{array}
\end{equation}
is a Brownian motion under $\P^u$. It means that $(x_t)_{t\le T}$ is a weak solution of the standard functional SDE \eqref{eds-contr}.

Next let $\t$ be a stopping time. We associate with the pair $(u,\t)$ a payoff $J(u,\t)$ given by :
\begin{equation}\label{pay}
J(u,\tau)=\E^u\left[\int_{0}^{\tau}\Gamma(s,x,u_s)ds+L_\tau\ind_{\{\tau<T\}}+g(x)\ind_{\{\tau=T\}}\right].
\end{equation}
The problem is to find a pair
$(u^*,\tau^*)$ which maximizes $J(u,\tau)$, i.e.,
\begin{equation}\lb{ineqoptim}
J(u^*,\tau^*)\ge J(u,\tau) \mx{ for any }(u,\t).
\end{equation}
This is the mixed control problem associated with $x$, $f$, $\Gamma$ and $g$. It combines control and stopping. One can think of an agent who controls a system by choosing the probability $\P^u$ which fixes its dynamics, weak solution of \eqref{eds-contr},  up to the time when she/he makes the decision to stop control at time $\t$. Therefore its payoff is given by $J(u,\t)$ and the problem is to find an optimal strategy $(u^*,\t^*)$, i.e., which satisfies \eqref{ineqoptim}.
\bigskip

Now let us introduce the Hamiltonian function $H$ associated with the mixed control problem which is given by: \begin{equation}\label{Hamiltonien}H(t,x,z,u):= z\sigma^{-1}(t,x)f(t,x,u)+\Gamma(t,x,u),\quad \forall(t,x,z,u)\in [0,T]\times\Omega\times\R^{d}\times
\ca.\end{equation} Note that for any $z\in \R^d$ and $u\in \cu$, the process $(H(t,x,z,u_t))_{t\leq
T}$ is ${\cal P}$-measurable. On the other hand thanks to Benes's selection Theorem (\cite{B}, Lemma 5, pp. 460), there exists a $\cp\otimes {\bf B}(\R^d)/{\bf B}(\ca)$-measurable function $u^*(t,x,z)$ such that for
any given $(t,x,z)\in[0,T]\times \Omega \times \R^d$,
\begin{equation}\label{BN}
H^*(t,x,z):=\sup\limits_{u\in \ca}H(t,x,z,u)=H(t,x,z,u^*(t,x,z)).\end{equation}

First we are going to show that the reflected BSDE associated with $(H^*,g,L)$ has a unique solution which moreover verifies some integrability properties. 
The following result is a step forward the proof of this result.
\begin{proposition}\label{BQ1m} For any $m\ge 1$, there exist $\cal P$-measurable processes $(Y^{*m},Z^{*m},K^{*m})$ in $\R^{1+d+1}$ such that:
\begin{enumerate}
\item[i)]For any $\gamma\geq 1$ and any stopping time $\tau\in[0.T]$,
\begin{equation}\label{m}
 \E\left[|Y^{*m}_\tau|^\gamma\right]\leq c,~~~~\forall m\geq 1,
\end{equation}
where $c$ is a constant independent of $m$ and $\tau$.
\item[ii)]
The triple
$(Y^{*m},Z^{*m},K^{*m})$ is a solution of the reflected BSDE associated with $(H^{*m},g,L)$ (in the sense of Def. \ref{defiLpsol}) where $$H^{*m}(s,x,z)={H^*}^+(s,x,z)-{H^*}^-(s,x,z)\ind_{\{1+||x||_s\leq m\}},$$ with for any 
$\a\in \R$, $\a^+=\a\vee 0$ and $\a^-=(-\a)\vee 0$.
\end{enumerate}
\end{proposition}
\ni \textbf{Proof}: For any $m,n\geq 1$ let us set,
\begin{equation}\label{Hnm}
H^{*n,m}(t,x,z)={H^*}^+(s,x,z)\ind_{\{1+||x||_s\leq n\}}-{H^*}^-(s,x,z)\ind_{\{1+||x||_s\leq m\}}.
\end{equation}
Then $H^{*n,m}$ is Lipschitz w.r.t $z$, therefore by El-Karoui et al's result \cite{EKPPQ}, there exists a triplet of processes $(Y^{*n,m},Z^{*n,m},K^{*n,m})\in \cs^2\times \H^2\times \cs^2 $ that satisfies:
\begin{equation}\label{BQ1nm}
\left\{\begin{array}{l} Y^{*n,m}_t=g(x)+\int_{t}^{T}H^{*n,m}(s,x,Z^{*n,m}_s)ds+K^{*n,m}_T-K^{*n,m}_t-\int_{t}^{T}Z^{*n,m}_s dB_s,\\
L_t \leq Y^{*n,m}_t, \quad \forall t\in[0,T]\mx{ and }\int_{0}^{T}(Y^{*n,m}_s-L_s)dK^{*n,m}_s=0.
\end{array}\right.
\end{equation}
Thus from the definition of $H^{*n,m}$, we can clearly see that it is a non-decreasing (resp. non-increasing) sequence of functions w.r.t $n$ (resp. $m$). Then once again by comparison (see \cite{EKPPQ}) we deduce that
\begin{equation}\label{inegynm}Y^{*n,m}\leq Y^{*n+1,m} \mbox{ and }Y^{*n,m+1}\leq Y^{*n,m}.\end{equation}
On the other hand we have
\begin{eqnarray*}
&& H^{*n,m}(s,x,Z^{*n,m}_s)\leq C[(1+||x||_s)\wedge n]|Z^{*n,m}_s|+C[(1+||x||^p_s).
\end{eqnarray*}
Then again by comparison we have: $\forall t\in[0,T]$
\begin{equation}\label{LY}
L_t\leq Y^{*n,m}_t\leq Y^n_t\leq Y_t,~~~~\forall n,m\geq 1,
\end{equation}
where $Y^n$ (resp. $Y$) is the process of (\ref{impn}) (resp. of Proposition \ref{unique}).\\

Next we will divide the proof into two steps. In the first (resp. second)  one we will prove i) (resp.  ii)).
\ms

\noindent \underline{Proof of i)}: From the growth condition on $L$ \eqref{grcdtl}, there exist  positive constant $C$ and $p$ such that $\forall t\in[0,T]$ and $x\in \Omega$
$$|L_t|\leq C(1+||x||^p_t).$$
From Lemma \ref{lemme1} we have that for any $\gamma\geq 1$ and any stopping times $\tau\in[0,T]$,
$$\E\left[|Y_\tau|^\gamma\right]\leq c, $$
where $c$ is a constant independent of $\tau$.
Then from (\ref{LY}), we deduce that
for any $\gamma\geq 1$ and any stopping time $\tau\in[0.T]$,
\begin{equation}\label{nm}
\E\left[|Y^{*n,m}_\tau|^\gamma\right]\leq c,~~~~\forall n,m\geq 1
\end{equation}
where $c$ is a constant that does not depend neither on $n$ nor $m$. Next for $m\ge 1$, let us set $$Y^{*m}=\liminf_{n\rightarrow +\infty}Y^{*n,m}.$$ Then by Fatou's Lemma and (\ref{nm}) we have:
\begin{equation}
\forall m\geq 1,\qquad \E\left[|Y^{*m}_\tau|^\gamma\right]\leq c \,\,\mbox{ and then }|Y^{*m}|<\infty,\,\, \P-a.s..
\end{equation}
On the other hand from \eqref{inegynm} it holds that
\begin{equation}\lb{compym}
L\le Y^{*,m+1}\le Y^{*,m}\le Y, \quad \fr m\ge 1.
\end{equation}
\underline{Proof of ii)}: Let $(\tau^{*}_k)_{k\geq 1}$ be the sequence of stopping times defined as follow:
$$\forall k\ge 1, \, \tau^{*}_k:=\inf\{t\geq 0, |Y_t|+\|x\|_t+|L_t|\ge \theta^{*}_k:=\theta^{*} +k \}\wedge T,$$
where $ \theta^{*}=|x_0|+|L_0|+|Y_0|$.\\
The sequence of stopping times  $(\tau^{*}_k)_{k\geq 1}$ is non-decreasing, of stationary type converging to $T$ since the processes $Y$, $x$ and $L$ are continuous. Moreover  for any $k\ge 1$,
$$ max\{\sup_{t\leq \tau^{*}_k}|L_t|,\sup_{t\leq \tau^{*}_k}|Y_t|,\sup_{t\leq \tau^{*}_k}|Y^{*n,m}_t|\} \leq  \theta^{*}_k$$\\
Let us now consider the following reflected BSDE. For any $n,m\geq 1$ and any $k\geq 1$ we have: $\forall t\in[0,T]$,
$$\left\{ \begin{array}{l}Y^{*n,m}_{t\wedge\tau^{*}_k} = Y^{*n,m}_{\tau^{*}_k} + \int_{t\wedge\tau^{*}_k}^{\tau^{*}_k}
 H^{*n,m}(s,x,Z^{*n,m}_s)ds +K^{*n,m}_{\tau^{*}_k}-K^{*n,m}_{t\wedge\tau^{*}_k}- \int_{t\wedge\tau^{*}_k}^{\tau^{*}_k} Z^{*n,m}_s dB_s;
\\\\
L_{t\wedge\tau^{*}_k}\leq Y^{*n,m}_{t\wedge\tau^{*}_k}\mbox{ and }\int_{0}^{\t^{*}_k}\left(Y^{*n,m}_s-L_s\right)dK^{*n,m}_s=0.
\end{array}\right.$$
Now we take into account ({\ref{estimx}}) and the fact that $\forall t\in[0,T]$, $|Y^{*n,m}_{t\wedge \tau^{*}_k}| \leq \theta^{*}_k$ and  $Y^{*0,m}\leq Y^{*n,m}\leq Y^{*m}$, then using It\^{o}'s formula with $\left(Y^{*n,m}_{t\wedge \tau^{*}_k}\right)^2$ to conclude that there exists a constant $C_k$ such that :
\begin{equation}\label{estimZk}\ba{l}
\E \left[\int_{0}^{\tau^{*}_k}|Z^{*n,m}_s|^2ds\right]\leq C_k.\ea
\end{equation}Next by It\^{o}'s formula we have: $\forall t\in[0,T]$
\begin{align}
(Y^{*n_1,m}_{t\wedge\tau^{*}_k}-Y^{*n_2,m}_{t\wedge\tau^{*}_k})^{2}&=(Y^{*n_1,m}_{\tau^{*}_k}-Y^{*n_2,m}_{\tau^{*}_k})^{2}-2\int_{t\wedge\tau^{*}_k}^{\tau^{*}_k}
(Y^{*n_1,m}_s-Y^{*n_2,m}_s)(Z^{*n_1,m}_s-Z^{*n_2,m}_s)dB_s\nn\\ \nonumber
& +2\int_{t\wedge\tau^{*}_k}^{\tau^{*}_k}(Y^{*n_1,m}_s-Y^{*n_2,m}_s)(H^{*n_1,m}(s,x,Z_s^{*n_1,m})-H^{*n_2,m}(s,x,Z_s^{*n_2,m}))ds\nn\\ \nonumber
&  +2\int_{t\wedge\tau^{*}_k}^{\tau^{*}_k}(Y^{*n_1,m}_s-Y^{*n_2,m}_s)d(K^{*n_1,m}_s-K^{*n_2,m}_s)-\int_{t\wedge\tau^{*}_k}^{\tau^{*}_k}|Z^{*n_1,m}_s-Z^{*n_2,m}_s|^2
ds.
\end{align}
Obviously for any $t\le T$, $$\ba{l}\int_{t\wedge\tau^{*}_k}^{\tau^{*}_k}(Y^{*n_1,m}_s-Y^{*n_2,m}_s)d(K^{*n_1,m}_s-K^{*n_2,m}_s)\leq 0.\ea$$
On the other hand, $$|H^{*n,m}(s,x,Z^{*n,m}_s)|\leq |H^*(s,x,Z^{*n,m}_s)|\leq C(1+||x||_s)|Z^{*n,m}_s|+C(1+||x||_s^p).$$
Therefore,
\begin{eqnarray}\label{conv}
&& (Y^{*n_1,m}_{t\wedge\tau^{*}_k}-Y^{*n_2,m}_{t\wedge\tau^{*}_k})^{2}+\int_{t\wedge\tau^{*}_k}^{\tau^{*m}_k}|Z^{*n_1,m}_s-Z^{*n_2,m}_s|^2
ds\leq \\ \nonumber
&&\qquad\qquad(Y^{*n_1,m}_{\tau^{*}_k}-Y^{*n_2,m}_{\tau^{*}_k})^{2}-2\int_{t\wedge\tau^{*}_k}^{\tau^{*}_k}(Y^{*n_1,m}_s-Y^{*n_2,m}_s)(Z^{*n_1,m}_s-Z^{*n_2,m}_s)dB_s\\ \nonumber
&& \qquad\qquad+2\int_{t\wedge\tau^{*}_k}^{\tau^{*}_k}|Y^{*n_1,m}_s-Y^{*n_2,m}_s|\left[C(1+||x||_s)\left(|Z^{*n_1,m}_s|+|Z^{*n_2,m}_s|\right)+2C(1+\|x\|^p_s)\right]ds.
\end{eqnarray}
Take now expectation in \eqref{conv} to deduce that:
\begin{multline*}
\E\left[\int_{t\wedge\tau^{*}_k}^{\tau^{*}_k}|Z^{*n_1,m}_s-Z^{*n_2,m}_s|^2 ds\right]\leq \E\left[(Y^{*n_1,m}_{\tau^{*}_k}-Y^{*n_2,m}_{\tau^{*}_k})^{2}\right]\\ \qquad\qquad
+2\E\left[\int_{t\wedge\tau^{*}_k}^{\tau^{*}_k}|Y^{*n_1,m}_s-Y^{*n_2,m}_s|\left[C(1+||x||_s)\left(|Z^{*n_1,m}_s|+|Z^{*n_2,m}_s|\right)+2C(1+\|x\|^p_s)\right]ds\right].
\end{multline*}
Then we conclude that
\begin{equation}\label{convz}\ba{l}
\E\left[\int_{t\wedge\tau^{*}_k}^{\tau^{*}_k}|Z^{*n_1,m}_s-Z^{*n_2,m}_s|^2 ds\right]\rightarrow 0\qquad as\qquad n_1,n_2\rightarrow +\infty.\ea
\end{equation}
This is due to estimate \eqref{estimZk}, the definition of $\tau^{*}_k$ and Cauchy-Schwarz inequality. Thus there exists a process $Z^{*m}_k\in\H^{2,d}$ such that $\left((Z^{*n,m}_t\ind_{\{0\leq t\leq \tau^{*}_k\}})_{t\leq T}\right)_{n\geq 1}\rightarrow (Z^{*m}_k(t))_{t\leq T}$ in $\H^{2,d}$ as $n\rightarrow +\infty$, $\forall\, m\geq 1$. Next by \eqref{conv} and the use of BDG inequality we obtain $$\E\left[\sup_{s\leq T} |Y^{*n_1,m}_{s\wedge\tau^{*}_k}-Y^{*n_2,m}_{s\wedge\tau^{*}_k}|^{2}\right] \rightarrow 0\qquad as\qquad n_1,n_2\rightarrow +\infty.$$
Since $(Y^{*n,m})_{n \geq 0}\rightarrow Y^{*m}$, the process $(Y^{*m}_{t\wedge\tau^{*}_k})_{t\geq 0}$ is continuous for every $k$ and $m$, also
\begin{equation}\label{convY}
\E\left[\sup_{s\leq T} |Y^{*n,m}_{s\wedge\tau^{*}_k}-Y^{*m}_{s\wedge\tau^{*}_k}|^{2}\right] \rightarrow 0\qquad as \qquad n\rightarrow +\infty.
\end{equation}
Next due to the fact that $(\tau^{*}_k)$ is of stationary type, then the process $Y^{*m}$ is continuous.\\
Now for $m\ge 1$, let us set
$$
K^{*m}_k(t)=Y^{*m}_0-Y^{*m}_{t\wedge \tau^{*}_k}- \int_0^{t\wedge\tau^{*}_k}H^{*m}(s,x,Z^{*m}_k(s))ds +\int^{t\wedge\tau^{*}_k}_0 Z^{*m}_k(s) dB_s.
$$
In taking into account \eqref{estimZk} and \eqref{convz} we deduce that for any $k\geq 1$, $m\geq 1$,\\   $\{(H^{*n,m}(t,x,Z^{*n,m}_t)\ind_{\{t\leq \tau^{*}_k\}})_{t\leq T}\}_{n\geq 1} $ converges in $\H^{2,d} $ to $(H^{*m}(t,x,Z^{*m}_t)\ind_{\{t\leq \tau^{*}_k\}})_{t\leq T},$
and then,
\begin{equation}\label{convKk}
\E\{\sup_{s\leq T}|K^{*n,m}_{s\wedge\tau^{*}_k}-K^{*m}_k(s)|^2\}\rightarrow 0\qquad as \qquad n\rightarrow +\infty.
\end{equation}
The uniform convergences in  \eqref{convY} and \eqref{convKk} imply that
$$\bal\int_{0}^{\tau^{*}_k}\left(Y^{*m}_s-L_s\right)dK^{*m}_k(s)=0.\ea$$
Now if we consider the reflected BSDE satisfied by $(Y^{*m},Z^{*m}_k,K^{*m}_k)$ on $[0,\tau^{*}_k]$ and the reflected BSDE satisfied by $(Y^{*m},Z^{*m}_{k+1},K^{*m}_{k+1})$ on $[0,\tau^{*}_k]$, we get by uniqueness: $\forall k\geq 1$, $\forall m\geq 1$
$$Z^{*m}_k(s)\ind_{\{s\leq \tau^{*}_k\}}=Z^{*m}_{k+1}(s)\ind_{\{s\leq \tau^{*}_k\}}dt\otimes d\P-a.e\mbox{ and }K^{*m}_{k+1}(s\wedge\tau^{*}_k)=K^{*m}_{k}(s\wedge\tau^{*}_k),\forall s\leq T.$$
So let us define $Z^{*m}$ and $K^{*m}$ by concatenation as follows: $\forall t\leq T$ and $m\geq 1$
\begin{equation}\label{defZm}
Z^{*m}_t=Z^{*m}_1(t)1_{\{t\le \t^{*}_1\}}+\sum_{k\ge 2}Z^{*m}_k(t)1_{\{\t^{*}_{k-1}< t\le \t^{*}_{k}\}}
\end{equation}
and
\begin{equation}\label{defKm}
K^{*m}_t=\left\{\begin{array}{l}
K^{*m}_1(t) \mbox{ if }t\leq \t^{*}_1;\\\\
(K^{*m}_{k+1}(t)-K^{*m}_{k+1}(\t^{*}_k)) +K^{*m}(\t^{*}_k)\mbox{ for }\t^{*}_k<t\leq \t^{*}_{k+1}, \,\,k\geq 1.
\end{array}\right.
\end{equation}
Then for any $k\geq 1$, and  any $m\geq 1$: $\frt$,
\begin{equation}\label{eqk12}\left\{ \begin{array}{l}Y^{*m}_{t\wedge\tau^{*}_k} = Y^{*m}_{\tau^{*}_k} + \int_{t\wedge\tau^{*}_k}^{\tau^{*}_k}
 H^{*m}(s,x,Z^{*m}_s)ds +K^{*m}_{\tau^{*}_k}-K^{*m}_{t\wedge\tau^{*}_k}- \int_{t\wedge\tau^{*}_k}^{\tau^{*}_k} Z^{*m}_sdB_s;
\\\\
L_{t\wedge\tau_k}\leq Y^{*m}_{t\wedge\tau_k}\mbox{ and }\int_{0}^{\t_k}\left(Y^{*m}_s-L_s\right)dK^{*m}_s=0.
\end{array}\right.\end{equation}
It is worth noticing that the processes $Z^{*m}$ and $K^{*m}$ are well defined since the sequence $(\tau^{*}_k)_{k\geq 1}$ is of stationary type. Moreover $K^{*m}$ is continuous, non-decreasing and $\P$-a.s. $K^{*m}_T(\omega)<+\infty$ and $\P$-a.s. $(Z^{*m}(t))_{t\leq T}$ is $dt$-square integrable. Finally going back to equation \eqref{eqk12}, take the limit when $k\rightarrow+\infty$ to obtain: $\forall m\geq 1$
\begin{equation}\label{eq12}\left\{ \begin{array}{l}Y^{*m}_{t} = g(x) + \int_{t}^{T}
 H^{*m}(s,x,Z^{*m}_s)ds +K^{*m}_{T}-K^{*m}_{t}- \int_{t}^{T} Z^{*m}_sdB_s;
\\\\
L\leq Y^{*m}\mbox{ and }\int_{0}^{T}\left(Y^{*m}_s-L_s\right)dK^{*m}_s=0.
\end{array}\right.\end{equation}
The proof is now complete.\qed
\ms

In the following result we show that the reflected BSDE \eqref{EK1intro} has a unique solution. Then after, we address the question of the link between the component $Y^*$ of the solution and the value function of the control-stopping problem. 
\begin{theorem}\label{BQ1}There exist $\cal P$-measurable processes $(Y^{*},Z^{*},K^{*})$ in $\R^{1+d+1}$ such that:

\item[i)] $Y^{*} \mx{ is continuous},
K^{*} \mx{ is continuous non-decreasing }(K^{*}_0=0) \mx{ and }
\P-a.s., \,\,(Z^{*m}(t))_{t\leq T}  \mx{ is }$\\$dt$-square 
integrable.

\item[ii)] For any $t\in[0,T]$,
\begin{equation}\label{BQm}
\left\{ \begin{array}{l}
Y^{*}_t=g(x)+\int_{t}^{T}H^{*}(s,x,Z^{*}_s)ds+K^{*}_T-K^{*}_t-\int_{t}^{T}Z^{*}_s dB_s,\\
L_t\leq Y^{*}_t ~and~\int_{0}^{T} (Y^{*}_s-L_s)dK^{*}_s=0.
\end{array}\right.
\end{equation}
\item[iii)] For any $\gamma\geq 1$ and any stopping time $\tau\in [0,T]$,
\begin{equation}\label{majy}
 \E\left[|Y^{*}_\tau|^\gamma\right]\leq c,
\end{equation}
where $c$ which does not depend on $\tau$.
\item[iv)]If $(\bar Y,\bar Z,\bar K)$ is another triple which satisfies i), ii) and iii), then
$(\bar Y,\bar Z,\bar K)=(Y^{*},Z^{*},K^{*})$, i.e., the solution
of the reflected BSDE associated with
$(g,H^*, L)$ is unique to satisfy i)-iii).
\end{theorem}
\ni \bf{Proof}: From Proposition \ref{BQ1m}, we have that for any $m\geq 1$, there exists a triplet $(Y^{*m},Z^{*m},K^{*m})$ that satisfies the following reflected BSDE: $\frt$,
$$\left\{\begin{array}{l} Y^{*m}_t=g(x)+\int_{t}^{T}H^{*m}(t,x,Z^{*m}_s)ds+K^{*m}_T-K^{*m}_t-\int_{t}^{T}Z^{*m}_s dB_s,\\
L_t \leq Y^{*m}_t, \mx{ and }\int_{0}^{T}(Y^{*m}_s-L_s)dK^{*m}_s=0.
\end{array}\right.$$
where $H^{*m}(s,x,z)=H^{*+}(s,x,z)-H^{*-}(s,x,z)\ind_{\{1+||x||_s\leq m\}}$.\\
From \eqref{compym} we know that for any $m\ge 1$,
$
L\le Y^{*,m+1}\le Y^{*,m}\le Y$. So, let us set for $t\le T$,
$$Y^*_t=\lim_m Y^{*,m}_t.$$
Therefore $L\le Y^*\le Y$ and then by  \eqref{estimx}, \eqref{grcdtl} and \eqref{ygamma}, $Y^*$ verifies \eqref{majy}.

Next let $(\tau^*_k)_{k\geq 1}$ be the sequence of stopping times defined as follows:
$$ \forall k\ge 1,~~~\tau^{*}_k:=\inf\{t\geq 0, |Y_t|+\|x\|_t+|L_t|\ge \theta^{*} +k \}\wedge T,$$
where $ \theta^{*}=|x_0|+|L_0|+|Y_0|$.
Thanks to the continuity of $Y$, $L$ and $x$, the sequence of stopping times $(\tau^*_k)_{k\geq 1}$ is increasing of stationary type and converges to $T$.
In the same way as in the proof of Proposition \ref{BQ1m}, there exists a constant $C_k$ (depending on $k$) such that:
$$\bal \E [\int_{0}^{\tau^*_k}|Z^{*m}_s|^2ds] \leq C_k.\ea$$
This inequality follows in a classic way after using It\^{o}'s formula with $\left(Y^{*m}_{t\wedge\tau^*_k}\right)^2$   and taking into account \eqref{estimx}.
Next we apply again It\^{o}'s formula and we have:
\begin{align}
(Y^{*m}_{t\wedge\tau^{*}_k}-Y^{*n}_{t\wedge\tau^{*}_k})^{2}&=(Y^{*m}_{\tau^{*}_k}-Y^{*n}_{\tau^{*}_k})^{2}-2\int_{t\wedge\tau^{*}_k}^{\tau^{*}_k}(Y^{*m}_s-Y^{*n}_s)(Z^{*m}_s-Z^{*n}_s)dB_s\nn\\ \nonumber
& +2\int_{t\wedge\tau^{*}_k}^{\tau^{*}_k}(Y^{*m}_s-Y^{*n}_s)(H^{*m}(s,x,Z_s^{*m})-H^{*n}(s,x,Z_s^{*n}))ds\nn\\ \nonumber
&  +2\int_{t\wedge\tau^{*}_k}^{\tau^{*}_k}(Y^{*m}_s-Y^{*n}_s)d(K^{*m}_s-K^{*n}_s)-\int_{t\wedge\tau^{*}_k}^{\tau^{*}_k}|Z^{*m}_s-Z^{*n}_s|^2
ds.
\end{align}
It follows that (in the same way as in the proof of Proposition \ref{BQ1m}): for any $k\ge 1$,
\begin{enumerate}
\item[a)] There exists a process $Z^*_k\in\H^{2,d}$, such that $((Z^{*m}_s\ind_{\{0\leq s\leq \tau^*_k\}})_{s\leq T})_{m\geq 1}\rightarrow_mZ^*_k(s)$ in $\H^{2,d}$, which is a consequence of the fact that  $((Z^{*m}\ind_{\{0\leq s\leq \tau^*_k\}})_{s\leq T})_{m\geq 1} $ is a Cauchy sequence in $\H^{2,d}$.
\item[b)]
$\lim_{m\rightarrow +\infty}\E\left[\sup_{s\leq T}|Y^{*m}_{s\wedge\tau^*_k}-Y^*_{s\wedge\tau^*_k}|^2\right]=0$ and the process $Y^*$ is continuous.
\end{enumerate}
Next we set
$$K^{*}_{k}(t)=Y^{*}_0-Y^{*}_{t\wedge\tau^*_k}-\int_{0}^{t\wedge{\tau^*_k}}
H^*(s,x,Z^{*}_k(s))ds+\int_{0}^{t\wedge\tau^*_k}Z^{*}_k(s)dB_s,t\le T,$$
and then we have,
\begin{equation}\label{convK*}
\lim_{m\rightarrow +\infty}\E\left[\sup_{s\leq T}|K^{*m}_{s\wedge\tau^*_k}-K^{*}_{k}(s)|^2\right]=0.
\end{equation}
It follows, from the uniform convergence of $(Y^{*m})_{m\geq 1}$ and $(K^{*m})_{m\geq 1}$ (see \cite{KF70}, pp. 370), that:
$$\bal \int_{0}^{\tau^*_k}(Y^*_s-L_s)dK_s=0\ea$$
Now considering the reflected BSDE satisfied by $(Y^*,Z^*_k,K^*_k)$ and the one satisfied by $(Y^*,Z^*_{k+1},K^*_{k+1})$ on $[0,\t^*_k]$, yields. For any $t\leq T$,
\begin{equation}\label{eqk1*}\left\{ \begin{array}{l}Y^*_{t\wedge\tau^*_k} = Y_{\tau^*_k} + \int_{t\wedge\tau^*_k}^{\tau^*_k}
 H^*(s,x,Z^*_k(s))ds +K^*_k({\tau^*_k})-K^*_k({t\wedge\tau^*_k})- \int_{t\wedge\tau^*_k}^{\tau^*_k} Z^*_k(s)dB_s;
\\\\
L_{t\wedge\tau^*_k}\leq Y^*_{t\wedge\tau^*_k}\mbox{ and }\int_{0}^{\t_k}\left(Y^*_s-L_s\right)dK^*_k(s)=0
\end{array}\right.\end{equation}
and
\begin{equation}\label{eqk2*}\left\{ \begin{array}{l}Y^*_{t\wedge\tau^*_k} = Y^*_{\tau^*_k} + \int_{t\wedge\tau^*_k}^{\tau^*_{k}}
 H^*(s,x,Z^*_{k+1}(s))ds +K^*_{k+1}({\tau^*_k})-K^*_{k+1}({t\wedge\tau^*_k})- \int_{t\wedge\tau^*_k}^{\tau^*_k} Z^*_{k+1}(s)dB_s;
\\\\
L_{t\wedge\tau^*_k}\leq Y^*_{t\wedge\tau^*_k}\mbox{ and }\int_{0}^{\t^*_k}\left(Y^*_s-L_s\right)dK^*_{k+1}(s)=0.
\end{array}\right.\end{equation}
Therefore for any $k\ge 1$, by uniqueness, we have:
$$
Z^*_k(s)1_{\{s\le \t^*_k\}}=Z^*_{k+1}(s)1_{\{s\le \t^*_k\}} \,\,dt\otimes d\P-a.e \mbox{ and }K^*_{k+1}({s\wedge\tau^*_k})=K^*_{k}({s\wedge\tau^*_k}),\,\,\fr s\le T.
$$
Finally by concatenation let us define the processes $Z^*$ and $K^*$ as follows: $\forall t\le T$,
$$
Z^*_t=Z^*_1(t)1_{\{t\le \t^*_1\}}+\sum_{k\ge 2}Z^*_k(t)1_{\{\t^*_{k-1}< t\le \t^*_{k}\}}$$ and $$K^*_t=\left\{\begin{array}{l}
K^*_1(t) \mbox{ if }t\leq \t^*_1;\\\\
(K^*_{k+1}(t)-K^*_{k+1}(\t^*_k)) +K^*(\t^*_k)\mbox{ for }\t^*_k<t\leq \t^*_{k+1}, \,\,k\geq 1.
\end{array}\right.
$$
The processes $Z^*$ and $K^*$ are well-defined. This is due to the fact that the sequence ($\t^*_k)_{k\ge 1}$ is of stationary type. Moreover $K^*$ is continuous non-decreasing and $\P$-a.s. $K^*_T(\o)<+\infty$ and $(Z^*_s(\o))_{s\le T}$ is $ds$-square integrable. Finally for any $k\ge 1$ it holds:
\begin{equation}\lb{eqk1quatre}\left\{ \begin{array}{l}Y^*_{t\wedge\tau^*_k} = Y_{\tau^*_k} + \int_{t\wedge\tau^*_k}^{\tau^*_k}
 H^*(s,x,Z^*_s)ds +K^*_{\tau^*_k}-K^*_{t\wedge\tau^*_k}- \int_{t\wedge\tau^*_k}^{\tau^*_k} Z^*_sdB_s;
\\\\
L_{t\wedge\tau^*_k}\leq Y^*_{t\wedge\tau^*_k}\mbox{ and }\int_{0}^{\t^*_k}\left(Y^*_s-L_s\right)dK^*_s=0.
\end{array}\right.\end{equation}
We now take the limit when $k\rightarrow+\infty$ and since once more ($\t^*_k)_{k\ge 1}$ is of stationary type, we have: $\fr t\le T$,
\begin{equation}\lb{eqk1cinq}\left\{ \begin{array}{l}Y^*_{t} = g(x)+ \int_{t}^{T}
 H^*(s,x,Z^*_s)ds +K^*_T-K^*_t- \int_{t}^{T} Z^*_s dB_s;
\\\\
L_{t}\leq Y^*_{t}\mbox{ and }\int_{0}^{T}\left(Y^*_s-L_s\right)dK^*_s=0,
\end{array}\right.\end{equation} which completes the proof of  $ii)$.
\ms

We will now prove $iv)$.
Let $(\bar{Y}, \bar{Z}, \bar{K})$ be another triple which satisfies i), ii) and iii). Then, using It\^{o}'s formula, we obtain:

\begin{equation} \label{Ystar}
\begin{array}{lll}
(Y^{*}_{t}-\bar{Y}_{t})^{2}&=-2\int_{t}^{T}(Y^*_s-\bar{Y}_s)(Z^*_s-\bar{Z}_s)dB_s\\
& \qquad +2\int_{t}^{T}(Y^{*}_s-\bar{Y}_s)(H^{*}(s,x,Z_s^{*})-H^{*}(s,x,\bar{Z}_s))ds \\
&  \qquad +2\int_{t}^{T}(Y^{*}_s-\bar{Y}_s)d(K^{*}_s-\bar{K}_s)-\int_{t}^{T}|Z^{*}_s-\bar{Z}_s|^2 ds.
\end{array}
\end{equation}
Next let $\P^*$ be the probability, equivalent to $\P$, defined as follows: $$d\P^*=L_T d\P$$ with
$$
L_T:=e^{\int_0^T \frac{H^{*}(s,x,Z_s^{*})-H^{*}(s,x,\bar{Z}_s)}{Z^*_s - \bar{Z}_s}1_{\{Z^*_s - \bar{Z}_s\neq 0\}}dB_s-\frac{1}{2}\int_0^T\|\frac{H^{*}(s,x,Z_s^{*})-H^{*}(s,x,\bar{Z}_s)}{Z^*_s - \bar{Z}_s}1_{\{Z^*_s - \bar{Z}_s\neq 0\}}\|^2ds }
$$and where 
$$
\Delta H^*(s):=\frac{H^{*}(s,x,Z_s^{*})-H^{*}(s,x,\bar{Z}_s)}{Z^*_s - \bar{Z}_s}1_{\{Z^*_s - \bar{Z}_s\neq 0\}}
$$is a $\cp$-measurable, $\R^d$-valued stochastic process such that 
$$\fr s\leq T, \,\,
H^{*}(s,x,Z_s^{*})-H^{*}(s,x,\bar{Z}_s)=\Delta H^*(s)\times (Z_s^{*}-\bar{Z}_s).
$$
As 
$$|H^{*}(s,x,z)-H^{*}(s,x,\bar z)|\le C(1+\|x\|_s)|z-\bar z| 
$$
then 
$$\fr s\leq T, \,\,
|\Delta H^*(s)|\le C(1+\|x\|_s).
$$It means that $P^*$ is actually a probability equivalent to $\P$ (by Lemma \ref{lemmeproba}). Next for $k\ge 1$, let $\t_k$ be the following stopping time:
$$\begin{array}{l}\tau_k:=\inf\{t\geq 0, |Y^*_t|+|\bar{Y}_t|+\int_0^t|Z^*_s|ds+\int_0^t|\bar{Z}_s|ds \ge  k +
|Y^*_0|+|\bar{Y}_0|
\}\wedge T.\end{array}$$
As $(Y^*_s-\bar{Y}_s)d(K^*_s-\bar{K}_s)\le 0$ then going back to \eqref{Ystar} to obtain:
$$\bal (Y^{*}_{t\wedge \tau_k}-\bar{Y}_{t \wedge \tau_k})^{2}\leq (Y^{*}_{\tau_k}-\bar{Y}_{\tau_k})^{2} +2 \int_{t\wedge \tau_k}^{\tau_k}(\bar{Y}_s-Y^*_s)(Z^*_s-\bar{Z}_s)d\widetilde{B}_s\ea $$
where
$
(\widetilde{B}_t:=B_t-\int_0^t\Delta H^*(s)ds)_{t\le T}$
is a Brownian motion under $\P^*$. Thus for any $t\le T$ and $k\ge 1$,
\begin{equation}\lb{equalityybary}\E^{\P^*}[(Y^{*}_{t\wedge \tau_k}-\bar{Y}_{t \wedge \tau_k})^{2}]\leq \E^{\P^*}[(Y^{*}_{\tau_k}-\bar{Y}_{\tau_k})^{2}].\end{equation}
But for any $\gamma \ge 1$ and $\tau$ stopping time, $
\E[|Y_\tau|^{\gamma}+|Y^*_\tau|^{\gamma}]\leq C$ and by Lemma \ref{lemmeproba} there exists a constant $p>1$ such that $\E[L_T^p]<\infty$. Then
there exists a constant $C$ such that for any stopping time $\t$, $\E^{\P^*}[(Y^{*}_{\tau}-\bar{Y}_{\tau})^{2}]\leq C. $ Consequently, the process $(Y^{*}-\bar{Y})^2$ is of class [D] under the
probability $\P^*$. Therefore (one can see e.g. \cite{DMM1}, Theorem 21, pp. 36)
$$
\E^{\P^*}[(Y^{*}_{\tau_k}-\bar{Y}_{\tau_k})^{2}]\rightarrow_k \E^{\P^*}[(Y^{*}_{T}-\bar{Y}_{T})^{2}]=0.
$$
Going back now to \eqref{equalityybary}, using Fatou's Lemma to obtain
$\E^{\P^*}[(Y^{*}_{t}-\bar{Y}_{t})^{2}]=0$ for any $t\le T$. It implies that $Y^*=\bar Y$, $\P^*$ and $\P$-a.s. since the probabilities are equivalent. Thus we have also $Z=\bar Z$ and $K=\bar K$, i.e. uniqueness. \qed
\ms

We then have the following result:
\ms

\begin{theorem}\lb{theorem32}
Let $(Y^*,Z^*,K^*)$ be the solution of the reflected BSDE associated with $(H^*,g(x),L)$, $u^*:=(u^*(t,x,Z^*_t))_{t\le T}$  and
$\tau^*=\inf\{t\in[0,T],Y^*_t\leq L_t\}\wedge T$. Then,
\begin{equation}\label{opti}
Y^*_0=J(u^*,\tau^*)=\sup_{u\in {\cal U},\tau\geq 0}J(u,\tau)
\end{equation}
i.e., $(u^*,\tau^*)$ is an optimal strategy of the mixed control problem.
\end{theorem}
\noindent \bf{Proof}:
Recall that $(Y^*,Z^*,K^*)$  verifies: $\frt$,
\begin{equation}\label{BM}
\left\{\begin{array}{l}
Y^*_t=g(x)+\int_{t}^{T}H^*(s,x,Z^*_s)ds+K^*_{T}-K^*_t-\int_{t}^{T}Z^*_sdB_s,\\
L\leq Y^* \mx{ and }\int_{0}^{T} (Y^*_s-L_s)dK^*_s=0.
\end{array}\right.
\end{equation}Next for $k\geq 1$, let us set
$$\bal\gamma_k=\inf\{s\geq 0, |L_s|+|Y^*_s|+\int_0^s|Z^*_r|^2dr \geq k+|L_0|+|Y^*_0| \}\wedge T.\ea$$
Now since $Y^*_0 $ is a deterministic constant we have:
$$\bal Y^*_0=\E^{u^*}\left[Y^*_{\tau^*\wedge\gamma_k}+\int_{0}^{\tau^*\wedge\gamma_k}H^*(s,x,Z^*_s)ds+K^*_{\tau^*\wedge\gamma_k}-\int_{0}^{\tau^*\wedge\gamma_k}Z^*_sdB_s\right],\ea$$
$$\bal \qquad\qquad\quad=\E^{u^*}\left[Y^*_{\tau^*\wedge\gamma_k}+\int_{0}^{\tau^*\wedge\gamma_k}\Gamma(s,x,u^*(s,x,Z^*_s))ds+K^*_{\tau^*\wedge\gamma_k}
-\int_{0}^{\tau^*\wedge\gamma_k}Z^*_sdB^{u^*}_s\right].\ea$$
From the definition of $\tau^*$ and the properties of reflected BSDEs we know that
the process $K^*_{\tau^*}$ does not increase between $0$ and $\tau^*$ then $K^*_{\tau^*\wedge\gamma_k}=0$. On the
other hand, using the Burkholder-Davis-Gundy inequality and the assumptions on $f$ we deduce that $\{\int_{0}^{t}Z^*_sdB^{u^*}_s, t\in [0,\gamma_k]\}$ is  an $\P^{u^*}$-martingale. Then,
$$\bal Y^*_0=\E^{u^*}\left[\int_{0}^{\tau^*\wedge\gamma_k}\Gamma(s,x,u^*(s,x,Z^*_s))ds+L_{\tau^*\wedge\gamma_k}\ind_{\{\tau^*\wedge\gamma_k<T\}}+g(x)\ind_{\{\tau^*\wedge\tau_k=T\}}\right].\ea$$But the sequence of stopping times $(\gamma_k)_{k\geq 1}$
is increasing, of stationary type and converges to T, then by taking the limit when $k\rightarrow+\infty$ we deduce that,
\begin{equation}\label{OPTP}
Y^*_0=J(u^*,\tau^*).
\end{equation}
Now, let $u$ be an admissible control and $\tau$ be an arbitrary stopping time. Since $\P$ and
$\P^{u} $ are  equivalent  probabilities on $(\Omega,\F)$ we have:
\begin{eqnarray*}
&& Y^*_0=\E^{u}[Y^*_0]=\E^{u}\left[Y^*_{\tau\wedge\gamma_k}+\int_{0}^{\tau\wedge\gamma_k}H^*(s,x,Z^*_s)ds+K^*_{\tau\wedge\gamma_k}-\int_{0}^{\tau\wedge\gamma_k}Z^*_sdB_s\right],\\ \nonumber
&& \qquad\qquad\qquad=\E^{u}\left[Y^*_{\tau\wedge\gamma_k}+\int_{0}^{\tau\wedge\gamma_k}\Gamma(s,x,u_s)ds+K^*_{\tau\wedge\gamma_k}-\int_{0}^{\tau\wedge\gamma_k}Z^*_sdB^{u}_s\right]\\ \nonumber
&& \qquad\qquad\qquad\qquad\qquad\qquad+\E^u\left[\int_{0}^{\tau\wedge\gamma_k}\left(H^*(s,x,Z^*_s)-H(s,x,Z^*_s,u_s)\right)\right].
\end{eqnarray*}
But $K^*_{\tau\wedge\gamma_k}\geq 0$, $H^*(s,x,Z^*_s)-H(s,x,Z^*_s,u_s)\geq 0$ and $\{\int_{0}^{t}Z^*_sdB^{u}_s, t\in [0,{\gamma_k}]\}$ is a $\P^{u}$-martingale, then,
$$Y^*_0\geq \E^{u}\left[\int_{0}^{\tau\wedge\gamma_k}\Gamma(s,x,u_s)ds+L_{\tau\wedge\gamma_k}\ind_{\{\tau\wedge\tau_k<T\}}+g(x)\ind_{\{\tau\wedge\gamma_k=T\}}\right].$$
The sequence of stopping times $(\gamma_k)_{k\geq 1}$
is increasing, of stationary type and converges to T, then by taking the limit as $k\rightarrow+\infty$ we get,
$$Y^*_0=J(u^*,\tau^*)\geq\E^{u}\left[\int_{0}^{\tau}\Gamma(s,x,u_s)ds+L_{\tau}\ind_{\{\tau<T\}}+g(x)\ind_{\{\tau=T\}}\right]= J(u,\tau)$$which is the claim. \qed
\begin{remark}\lb{remarque31} The process $Y^*$ is the value function of the mixed optimal control problem, i.e., for any $t\in [0,T]$,  
\begin{equation}\lb{valprocess}
Y^*_t=\mx{esssup}_{{\tau\ge t, u\in \cu}}\E^{u}\left[\int_{t}^{\tau}\Gamma(s,x,u_s)ds+L_{\tau}\ind_{\{\tau<T\}}+g(x)\ind_{\{\tau=T\}}|\F_t\right].
\end{equation}This is another way to show uniqueness of the solution of \eqref{BQm}. \qed
\end{remark}
\section{The Markovian framework: The HJB equation associated with the mixed control problem}
In this section we are going to restrict the previous framework to the Markovian one and study the properties of the value function and the Hamilton-Jacobi-Bellman equation associated with the mixed control problem. We will show that the first one provides the unique solution in viscosity sense of the second. 

\ni To begin with let us introduce the following spaces:
\ms

\ni i) For any $\gamma \ge 0$, $\Pi_{\mbox{pg}}^\g$ be the following space of functions $\Psi$ defined on $\esp$, $\R$-valued and  such that $$\fr \txt , \,\,|\Psi(t,x)|\leq C(1+|x|^\gamma) \mbox{ for some constants } C\ge 0.$$

\ni ii) $\Pi_{\mbox{pg}}=\cup_{\g \ge 0}\Pi_{pl}^\g$.
\ms

Next let us specify the functions $\sigma$, $f$, $g$, $h$ and $\Gamma$, with their properties, in the following assumptions which we assume satisfied hereafter:
\ms

\ni \bf{(HM)}: the functions $\sigma$, $f$, $g$, $h$ and $\Gamma$ verify:
\\
\ms

\noindent a) $\sigma:[0,T]\times \R^d\rightarrow \R^{d\times d} $ is a bounded continuous function, Lipschitz w.r.t $x$ invertible and  its inverse $\sigma^{-1}(t,x)$ is bounded and continuous. 
\ms

\noindent b) $h: [0,T]\times \R^d\rightarrow \R$ and
$g: \R^d\rightarrow \R$. They are continuous and belong to $\pg$. 
\ms

\ni c) $\Gamma:[0,T]\times \R^d \times \ca\rightarrow \R$ is continuous. In addition $\sup_{a\in \ca}\Gamma$ belongs to $\pg$. 
\ms

\ni d) $f:[0,T]\times \R^d \times \ca\rightarrow \R^d$ is continuous and the mapping $x\in \R^d \mapsto f(t,x,a)$ is  Lipschitz uniformly w.r.t $(t,a)$. Moreover it is of linear growth, i.e., $\sup_{a\in \ca}f$ belongs to $\pg^1$ and then  for any $(t,x,a)$,
\begin{equation}\lb{fcl}
|f(t,x,a)|\le C(1+|x|) \mbox{ for some constant } C\ge 0. \qed
\end{equation}

Now for $\txt$, let $X^\tx:=(X^\tx_s)_{s\le T}$ be the solution of the following standard differential equation:
\begin{equation}\label{dyn}
\left\{
\begin{array}{l}
dX^{t,x}_s=\sigma(s,X^{t,x}_s)dB_s,~~s\in[t,T];\\
 X^{t,x}_s=x, \qquad\qquad\qquad\qquad s\leq t.
\end{array}\right.
\end{equation}In the
following result we collect some properties of $X^{t,x}$.
\begin{proposition}\label{estimxx} (see e.g. \cite{[RY]}) The process $X^{t,x}$ satisfies
the following estimates:
\begin{itemize}
\item [$(i)$] For any $q\geq 2$, there exists a constant $C$ such that
\begin{equation}\label{estimat1}
\E[\sup_{0\le s\leq T}|X^{t,x}_s|^q]\leq C(1+|x|^q).
\end{equation}
\item[$(ii)$] There exists a constant $C$ such that for any $t,t'\in [0,T]$ and $x,x'\in \R^d$,
\begin{equation}\label{estimat2}
\E[\sup_{0\le s\leq T}|X^{t,x}_s-X^{t',x'}_s|^2]\leq
C(1+|x|^2)(|x-x'|^2+|t-t'|). \Box
\end{equation}
\end{itemize}
\end{proposition}

Next for $(t,x,z)\in [0,T]\times \R^{d+d}$, let us introduce the function $H^*$ which is the same as the one given in \eqref{BN}  in this Markov setting:
\begin{equation}\label{BNmk}
H^*(t,x,z):=\sup\limits_{a\in \ca}H(t,x,z,a)\end{equation}
where (with $a\in \ca$), \begin{equation}\label{Hamiltonienmk}H(t,x,z,a):= z\sigma^{-1}(t,x)f(t,x,a)+\Gamma(t,x,a).\end{equation} 
The function $H^*$ is continuous in all its arguments (see Lemma \ref{conth*} below), property which is needed later to deal with the HJB equation associated with the mixed control problem.  Now to proceed, for $(t,x,z)\in [0,T]\times \R^{d+d}$ and $n,m\ge 0$, let us define ${\bar H}^{*n,m}$ by $$
\begin{array}{ll}{\bar H}^{*n,m}(t,x,z):={H^*}^+(t,x,z)\rho_n(x)-{H^*}^-(t,x,z)\rho_m(x).
\end{array}$$
where for an $x\in \R^d$, $$\rho_m(x)=
\ind_{\{|x| \leq m\}}+(x+1+m)\ind_{\{-m-1\leq x \leq -m\}}+(-x+1+m)
\ind_{\{m\leq x \leq m+1\}}.$$
The function ${\bar H}^{*n,m}$ is also a truncation of $H^*$ which is moreover continuous w.r.t $(t,x,z)$. It is uniformly Lipschitz in $z$ and have the same monotonicity properties as ${H}^{*n,m}$ defined in \eqref{Hnm}, i.e., ${\bar H}^{*n,m}$ is also increasing (resp. decreasing) w.r.t. $n$ (resp. $m$) and $\lim_{m\rw \infty}\lim_ {n\rw \infty}{\bar H}^{*n,m}=\lim_{n\rw \infty} \lim_{m\rw \infty}{\bar H}^{*n,m}=H^*$ since $\rho_m(x)\nearrow 1$. Finally the following estimate holds true: 
\begin{equation}\label{ineqhnm}
|{\bar H}^{*n,m}(t,x,z)|\leq \bar \varphi (t,x,z):=C(1+|x|)|z|+C(1+|x|^p),\,\,\forall (t,x,z)\in [0,T]\times \R^{d+d}.
\end{equation}
Next for $n,m\ge 0$, let  $(\bar{Y}^{t,x,n,m},\bar{Z}^{t,x,n,m},\bar{K}^{t,x,n,m})$ be the unique solution of the following BSDE associated with $({\bar H}^{*n,m},g, h)$: $\forall s\le T$,
\begin{equation}\label{visco}
\left\{
\begin{array}{l}
\bar{Y}^{t,x,n,m} , \bar{K}^{t,x,n,m}\in \cs^2 
\mbox { and }\bar{Z}^{t,x,n,m}\in \H^{2,d} \,\,(\bar{K}^{t,x,n,m} \mbox{ increasing and } \bar{K}^{t,x,n,m}_0=0) ;\\
\bar{Y}^{t,x,n,m}_s = g(X^{t,x}_T) + \int_s^T
{\bar H}^{*n,m}(u,X^{t,x}_u,\bar{Z}^{t,x,n,m}_u) du +\bar{K}^{t,x,n,m}_T-\bar{K}^{t,x,n,m}_s- \int_s^T \bar{Z}^{t,x,n,m}_u dB_u;\\\\
h(s,X^{t,x}_s)\leq \bar{Y}^{t,x,n,m}_s\mx { and }\int_{0}^{T} (\bar{Y}^{t,x,n,m}_u-h(u,X^{t,x}_u))d\bar{K}^{t,x,n,m}_u=0.
\end{array}\right.
\end{equation}
By comparison we have: $\forall n,m\ge 0$, 
\begin{equation}\label{ybarnm}
\bar{Y}^{t,x,n,m}\le \bar{Y}^{t,x,n+1,m}\le \bar{Y}^{t,x,n+1,m-1}.
\end{equation} On the other hand by Theorem 8.5 in \cite{EKPPQ}, since the framework is Markovian, then: 
\medskip

\nd i) $\bar{u}_{n,m}(t,x)=\bar{Y}^{t,x,n,m}_t$, $(t,x)\in[0,T]\times \R^d$, is a deterministic continuous function which is moreover solution in viscosity sense of the following parabolic PDE: 
\begin{equation}\lb{edpunm}\left\{\begin{array}{l}\min\left[\bar{u}_{n,m}(t,x)-h(t,x),-\partial_t \bar{u}_{n,m}(t,x)-\cl \bar{u}_{n,m}(t,x)-{\bar H}^{*n,m}(t,x,\sigma(t,x)\nabla_x\bar{u}_{n,m}(t,x))\right]=0, \\(t,x)\in[0,T)\times \R^d;\\\\
\bar{u}_{n,m}(T,x)=g(x),x\in \R^d.
\end{array}\right.\end{equation}
ii) For any $s\in [t,T]$, 
\begin{equation}\lb{repynm}
\bar{Y}^{t,x,n,m}_s=\bar{u}_{n,m}(s,X^{t,x}_s).
\end{equation}
Next for $m\ge 0$ let us set 
\begin{equation}\label{hbar*m}
\begin{array}{l}\bar{H}^{*m}(t,x,z):=
\lim_{n\rw \infty}\bar{H}^{*n,m}(t,x,z)={H^*}^+(t,x,z)-
{H^*}^-(t,x,z)\rho_m(x).
\end{array}\end{equation}
and 
\begin{equation}\lb{defym}
\bar Y^{t,x,m}=\lim_{n\rw \infty}\nearrow \bar{Y}^{t,x,n,m}.
\end{equation}
By Proposition \ref{BQ1m}, there exist a ${\cal P}$-measurable process $\bar Z^{t,x,m}$ and an increasing continuous process 
$\bar K^{t,x,m}$ such that $(\bar Y^{t,x,m},
\bar Z^{t,x,m},\bar K^{t,x,m})$ is a solution of the reflected BSDE associated with 
$(\bar {H}^{*m}(s,X^{t,x}_s,z),g(X^{t,x}_T),h(s,X^{t,x}_s))$, i.e., $\fr s\le T$,
$$\left\{\begin{array}{l} \bar Y^{t,x,m}_s=g(X^{t,x}_T)+\int_{s}^{T}\bar H^{*m}(r,X^{t,x}_r,\bar Z^{t,x,m}_r)dr+\bar K^{t,x,m}_T-\bar K^{t,x,m}_s-\int_{s}^{T}\bar Z^{t,x,m}_r dB_r;\\\\
h(s,X^{t,x}_s)\leq Y^{t,x,m}_s \mx{ and }\int_{0}^{T}(Y^{t,x,m}_r-h(r,X^{t,x}_r))dK^{t,x,m}_r=0.
\end{array}\right.$$
Note that $\bar Y^{t,x,m}$ verifies the estimate \eqref{m}, and $\P-a.s$, $\bar Z^{t,x,m}$ is $ds$-square integrable and $\bar K^{t,x,m}_T<\infty$. The inequalities \eqref{ybarnm} imply that  for any $\txt$, 
$$
\bar u_{n,m}(t,x)\le \bar u_{n+1,m}(t,x)\le \bar u_{n+1,m-1}(t,x).
$$Now for any $\txt$ and $m\ge 0$, let us set 
$$
\bar u^{m}(t,x)=\lim_{n\rw \infty}\bar u_{n,m}(t,x).$$ Then by \eqref{repynm} and \eqref{defym}, it holds that for any $\txt$ 
\begin{equation}\lb{repym}
\bar Y^{m,t,x}_s=\bar u^{m}(s,X^{t,x}_s),\,\,\fr s\in [t,T].
\end{equation}
Next once more as a consequence of \eqref{ybarnm}, the sequence $(\bar Y^{m,t,x})_{m\ge 0}$ is decreasing and then so is the sequence of deterministic functions $(\bar u^{m})_{m\ge 0}$. So let us define the process $Y^{t,x}$ and the deterministic function $u$ by:
$$Y^{t,x}=\lim_{m\rw \infty}\searrow \bar Y^{t,x,m}\mbox{ and }u(t,x)=\lim_{m\rw \infty}\searrow \bar u^{m}(t,x).$$
By \eqref{repym}, those latter objects are connected by the following relation:
\begin{equation}\lb{repyu}
\bar Y^{t,x}_s=u(s,X^{t,x}_s),\,\,\fr s\in [t,T].
\end{equation}As in Theorem \ref{BQ1}, there exists a pair of processes
$(Z^{t,x},K^{t,x})$ valued in $\R^{d+1}$ such that $(Y^{t,x},Z^{t,x},K^{t,x})$ is the unique solution of the following reflected BSDE: $\forall s\in[0,T]$,
\begin{equation}\label{solutionmarkov}
\left\{
\begin{array}{l}
Y^{t,x}_s = g(X^{t,x}_T) + \int_s^T
H^*(r,X^{t,x}_r,Z^{t,x}_r) dr +K^{t,x}_T-K^{t,x}_s- \int_s^T Z^{t,x}_r dB_r;\\\\h(s,X^{t,x}_s)\leq Y^{t,x}_s\mx{ and }\int_{0}^{T} (h(r,X^{t,x}_r)- Y^{t,x}_r)dK^{t,x}_r=0.
\end{array}\right.
\end{equation}
\begin{lemma} \lb{lemmecroissancepoly}There exist two positive constants $C$ and $p$ such that for any $\txt$ and $m\ge 0$:\\
i)
\begin{equation}\lb{croipolyumu}
|\bar u^{m}(t,x)|+|u(t,x)|\le C(1+|x|^p).
\end{equation}
ii)
$$
\E[\int_0^T\{|\bar Z^{t,x,m}_r|^2+|Z^{t,x}_r|^2\}dr]\le C(1+|x|^p).
$$
\end{lemma}
\noindent \bf{Proof:} i) Recall the function $\bar \varphi$ introduced above in \eqref{ineqhnm}:
$$
\bar \varphi (t,x,z):=C(1+|x|)|z|+C(1+|x|^p).
$$
By Proposition \ref{unique}, there exists a triplet of processes 
$(y^{t,x},z^{t,x},k^{t,x})$ such that:\\\\
a) $(y^{t,x}, k^{t,x})$ is $\R^{1+1}$-valued, 
$k^{t,x}$ is  non-decreasing and $k^{t,x}_0=0$,
$(z_s^{t,x}(\o))_{s\le T}$ is $\R^d$-valued and 
$ds$-square integrable $\P-a.s.$\\
b) For any constant $\gm\ge 1$ and $\t$ a stopping time valued in $[0,T]$,
\begin{equation}\lb{ygamma2}
\E[|y^{t,x}_\t|^\gm]<+\infty.
\end{equation}
c) For any  $s\le T$,
\begin{equation}\label{y2} \left\{
\begin{array}{l}
y^{t,x}_s= g(X^{t,x}_T) + \int_{s}^{T}\{
C|z^{t,x}_r|(1+|X^{t,x}_r|)+C(1+|X^{t,x}_r|^p)\}dr-\int_{s}^{T}z^{t,x}_rdB_r+k_T^{t,x}-k^{t,x}_s;\\ \\y^{t,x}_s\geq l^{t,x}_{s}:=
h( s,X_{
s}^{t,x})\mbox{ and } \,\,
\int_{0}^{T}(y^{t,x}_r-l^{t,x}_{r})dk^{t,x}_r=0.
\end{array}
\right. \end{equation}
Now by the standard comparison result of solutions of reflected BSDEs and thanks to \eqref{ineqhnm}, one has:
\begin{equation}\lb{inegalitex22}
l^{t,x}_{s}\le \bar Y^{t,x,n,m}_s\le y^{t,x}_s, \forall s\le T.
\end{equation}
Next let $\bar \P^{t,x}$  be the probability, equivalent to $\P$, defined as follows: $$d\bar\P^{t,x}=M^{t,x}_T d\P$$where
$$\bal
M^{t,x}_T:=\exp\{\int_0^TC\sigma^{-1}(r,X^{t,x}_r)\ell(z^{t,x}_r)C(1+|X^{t,x}_r|)dB_r-
\frac{1}{2}\int_0^T\|C\sigma^{-1}(r,X^{t,x}_r)\ell(z^{t,x}_r)(1+|X^{t,x}_r|)\|^2dr \}\ea
$$
where $\ell$ is bounded measurable function such that $\ell(z).z=|z|$, $\fr z=(z_i)_{i=1,\dots,d}\in \R^d$ (see the proof of Lemma \ref{lemme1} for its defintion). Under $\bar \P^{t,x}$, $X^\tx$ is a weak solution of the following SDE:
\begin{equation}\label{dyn2}
\left\{
\begin{array}{l}
dX^{t,x}_s=C\ell(z^{t,x}_s)C(1+|X^{t,x}_s|)ds+\sigma(s,X^{t,x}_s)d\bar B_s,~~s\in[t,T];\\
 X^{t,x}_s=x, \qquad\qquad\qquad\qquad s\leq t.
\end{array}\right.
\end{equation}
were $\bar B$ is a Brownian motion under $\bar \P^{t,x}$. Then it verifies the following estimate: $\fr q\ge 2$,
\begin{equation}\label{estimat1x}
{\bar \E}^{t,x}[\sup_{0\le s\leq T}|X^{t,x}_s|^q]\leq C(1+|x|^q)
\end{equation}
where ${\bar \E}^{t,x}$ is the expectation under ${\bar \P}^{t,x}$; the constant $C$ does not depend on $t,x$. Next writing the reflected BSDE\eqref{y2} under the probability $\bar \P^{t,x}$ reads:
For any  $s\le T$,
\begin{equation}\label{y21} \left\{
\begin{array}{l}
y^{t,x}_s= g(X^{t,x}_T) + \int_{s}^{T}C(1+|X^{t,x}_r|^p)dr-\int_{s}^{T}z^{t,x}_rd\bar B_r+k_T^{t,x}-k^{t,x}_s;\\ \\y^{t,x}_s\geq l^{t,x}_{s}:=
h( s,X_{
s}^{t,x})\mbox{ and } \,\,
\int_{0}^{T}(y^{t,x}_r-l^{t,x}_{r})dk^{t,x}_r=0.
\end{array}
\right. \end{equation}
Taking now into account of \eqref{ygamma}, we obtain as previously (see the proof of Theorem \ref{theorem32} or Remark \ref{remarque31}) the following representation for $y^{t,x}$: $\fr s\le T$,
$$\bal y^{t,x}_s=\esssup_{\tau \ge s}\bar \E^{t,x}[\int_s^\tau C(1+|X^{t,x}_r|^p)dr +l^{t,x}_\tau\ind_{\{\tau<T\}}+g(X^{t,x}_T)\ind_{\{\tau=T\}}|\F_s].\ea$$
Next by the polynomial growth of $h$ and $g$ we deduce that: $\fr s\le T$,
$$
|y^{t,x}_s|\le C (1+\bar \E^{t,x}[\sup_{s\le T}|X^{t,x}_s|^q|\F_s])
$$
for some constants $C$ and $q$. Then by \eqref{estimat1x}, we deduce that 
$$
\bar \E^{t,x}[|y^{t,x}_t|]\le C (1+|x|^p)
$$
for some fixed $C$ and $p$. Going back now to \eqref{inegalitex22}, take $s=t$, expectation w.r.t $\bar \P^{t,x}$ (which is equivalent to $\P$) and since $Y^{t,x,n,m}_t$ is deterministic one deduces that 
$$
h(t,x)\le \bar u_{n,m}(t,x)\le \bar \E^{t,x}[y^{t,x}_t]\le C(1+|x|^p).$$
The proof is now completed in taking the limit w.r.t $n$ then $m$ since $h$ is polynomial growth .

ii) It is obtained by the use of It\^o's formula with $|Y_s^{t,x}|^2$ and $|\bar Y_s^{m,t,x}|^2$ and in taking into account the representations \eqref{repym}  and \eqref{repyu}, the polynomial growths \eqref{croipolyumu} of $u$ and $\bar u^{m}$, estimate \eqref{estimat1} on $X^{t,x}$ and finally the fact that $Z^{t,x}_s=\bar Z^{m,t,x}_s=0$ for $s\in [0,t]$.\qed
\begin{remark}By the polynomial growth of $u$ and estimate \eqref{estimat1}, we have also 
\begin{equation}\label{estimytx}\E[\sup_{s\le T}|Y^{t,x}_s|^p]<\infty, \,\,\forall p \ge 1,\end{equation}\qed\end{remark}

Next let $\Phi(t,x,z)$ be a function from $[0,T]\times \R^{d+d}$ into $\R$ which we assume continuous in all its arguments. Let us now consider the following PDE with obstacle:
\begin{equation}\label{kkk}
\left\{
\begin{array}{l}
\min\left[v(t,x)-h(t,x),
-\frac{\partial v}{\partial t}(t,x)-\cl v(t,x)-\Phi(t,x,\nabla_xv(t,x)\sigma(t,x))\right]=0, \,\,(t,x)\in[0,T[ \times \R ^d;\\\\
v(T,x)=g(x),\, x\in\R^d,
\end{array}
\right.
\end{equation}
where $\nabla_x$ is the derivative w.r.t. $x$ and $\cl$ is second order partial differential operator associated with $X^{t,x}_.$ i.e
$$\cl=\frac{1}{2}\sum_{i,j=1,d}(\sigma\sigma^\top)_{ij}(t,x)\partial^2_{x_ix_j}.$$

In the case when we take $\Phi=H^*$, we obtain the HJB equation associated with the mixed control problem in the Markov framework.
\medskip

Next let us define the notion of a solution of \eqref{kkk}. Let $v$ be a funtion defined on $\esp$ which is moreover locally bounded on each $\txp$. We define the upper (resp. lower) semi-continuous enveloppe of $v$ by: $\fr \txt$,
$$
v^*(t,x):=\limsup_{(t',x')\rw (t,x), t'<T}v(t',x') \,\,(\mx{resp.}\,\,v_*(t,x):=\liminf_{(t',x')\rw (t,x), t'<T}v(t',x')).
$$

The definitions of the limiting parabolic superjet (resp. subjet) $\bar J^{2,+}v$ (resp. $\bar J^{2,-}v$) of an upper (resp. a lower) semi-continuous function $v$ defined on $\esp$ are given e.g. \cite{CIL}, pp.47, 11.

\begin{definition}Let $v$ be function defined on $\esp$, $\R$-valued and locally bounded: 
\ms

\ni a) It is said a viscosity supersolution (resp.
subsolution) of (\ref{kkk}) if:

(i) $v_*(T,x) \geq g(x)$ (resp. $v^*(T,x) \le g(x)$), $\fr x\in \R^d$;

(ii) For any $(t,x)\in
[0,T)\times \R^d$ and $(p,q,X)\in \bar J^{2,-} v_*(t,x)$ (resp. $\bar J^{2,+}
v^* (t,x)$),
$$\min \left\{v_*(t,x)-h(x),-p -\frac{1}{2}Tr[\sigma^\top X \sigma (t,x)] -\Phi(t,x,q\sigma(t,x))\right\}\geq 0 $$ (resp.
\begin{equation}\label{kkk1}min \left\{v^*(t,x)-h(x),-p -\frac{1}{2}Tr[\sigma^\top X \sigma (t,x)] -\Phi(t,x,q\sigma(t,x))\right\}\leq 0).\end{equation} 
b) It is
called a viscosity solution if it is both a viscosity subsolution and
supersolution.
\end{definition}
\begin{remark}There is another definition of the notion of viscosity solution which uses test functions which we will use sometimes later on (one can see e.g. \cite{CIL} in pp. 10-11) for more details. \qed  
\end{remark}

To begin with we are going to deal with the issue of comparison principle, and consequently uniqueness of the solution, for the PDE with obstacle \eqref{kkk}. For that let us introduce the following assumptions on the function $\Phi$:
\medskip

\ni\bf{\underline{(H$\Phi$)}}: 
\medskip

\ni i) 
For any $(t,x,z,z')\in [0,T]\times \R^{d+d+d}$,
\begin{equation}\label{cd1gfi}
 \Phi(t,x,z)-\Phi(t,x,z')\ge -C(1+|x|)|z-z'|.
\end{equation} 
ii)
For any $\kappa>0$, there exists a function $\Psi_\kappa$ from $[-2\kappa,2\kappa]$ into $\R^+$, continuous, $\Psi_\kappa(0)=0$ and such that for any $(t,z,z')\in [0,T]\times \R^{2d }$, $|x|\le \kappa$, $|y|\le \kappa$,
\begin{equation}\label{cd2gfi}
 \Phi(t,x,z\sigma(t,x))-\Phi(t,y,z'\sigma(t,y))\le C_\kappa(|z|+|z'|)|x-y|+C_\kappa|z-z'|+\Psi_\kappa(|x-y|)
\end{equation} 
where $C_\kappa$ is a positive constant which may depend on $\kappa$. \qed
\ms

As a preliminary result we have:  
\begin{lemma}\label{lemma-supersol}Assume that the function $\Phi$ verifies \bf{(H$\Phi$)}-i). If $v$ is a supersolution of (\ref{kkk}) which belongs to $\pg$, i.e.,   
$$
\fr \txt, |v(t,x)|\leq C(1+|x|^{2\gamma})
$$
for some constants $C$ and $\gamma$ non negative. Then there exists $\lambda_0>0$ such that for any $\lambda \geq \lambda_0$ and $\theta >0$,  $(v(t,x)+\theta e^{-\lambda t}(1+\mid x\mid^{2\gamma+2})$ is a supersolution for (\ref{kkk}) .
\end{lemma}
\noindent $Proof. $ We assume w.l.o.g. that the function
$v(t,x)$ is lsc. First note that the condition at $T$ holds true since 
$\theta e^{-\lambda T}(1+\mid x\mid^{2\gamma+2})\ge 0$. Next let $t<T$ and $\varphi \in \mathcal{C}^{1,2}$ be such that the
function $\varphi-(v+\theta e^{-\lambda t}(1+\mid
x\mid^{2\gamma+2})$ has a local maximum in $(t,x)$ which is equal to
0. Since $v(t,x)$ is a supersolution for
(\ref{kkk}), then we have:
\begin{eqnarray}
&&\min\bigg\{ v(t,x)-h(t,x),\nonumber\\&&\hspace{10mm}-\partial_{t}\Big(\varphi(t,x)-\theta e^{-\lambda t}(1+\mid x\mid^{2\gamma+2})\Big)-\frac{1}{2}Tr\Big[\sigma.\sigma^\top(t,x)D^{2}_{xx}\Big(\varphi(t,x)-\theta
e^{-\lambda t}\mid x\mid^{2\gamma+2}\Big)\Big]\nonumber\\&&\hspace{10mm}-\Phi(t,x,\nabla_x(\varphi(t,x)-\theta
e^{-\lambda t}\mid x\mid^{2\gamma+2})\sigma(t,x))\bigg\}
\geq 0.\nonumber
\end{eqnarray}
Hence
\begin{eqnarray}\label{supersol1}
(v(t,x)+\theta e^{-\lambda t}(1+\mid
x\mid^{2\gamma+2}))-h(t,x) \geq v(t,x)-h(t,x) \geq 0.
\end{eqnarray}
On the other hand:
\begin{equation*}
\begin{array}{l}
\left.-\partial_{t}\Big(\varphi(t,x)-\theta e^{-\lambda t}(1+\mid x\mid^{2\gamma+2})\Big)-\frac{1}{2}Tr\Big[\sigma.\sigma^\top(t,x)D^{2}_{xx}\Big(\varphi(t,x)-\theta e^{-\lambda t}\mid x\mid^{2\gamma+2}\Big)\Big] \right.\\
\left.-\Phi(t,x,\nabla_x(\varphi(t,x)-\theta
e^{-\lambda t}\mid x\mid^{2\gamma+2})\sigma(t,x))\right. \geq 0.
\end{array}
\end{equation*}
Therefore
\begin{align}
-\partial_{t}\varphi(t,x)&-\frac{1}{2}Tr\Big[\sigma.\sigma^\top(t,x)D^{2}_{xx}\varphi(t,x)\Big]-\Phi(t,x,\nabla_x(\varphi(t,x))\sigma(t,x)) \nonumber\\
&\geq \theta \lambda e^{-\lambda t}(1+\mid
x\mid^{2\gamma+2})-\frac{1}{2}\theta e^{-\lambda t}
Tr\Big[\sigma.\sigma^\top(t,x)D^{2}_{xx}\mid x\mid^{2\gamma+2}\Big] \nonumber\\
&\qquad + \Phi(t,x,\nabla_x(\varphi(t,x)-\theta
e^{-\lambda t}\mid x\mid^{2\gamma+2})\sigma(t,x))-\Phi(t,x,\nabla_x(\varphi(t,x))\sigma(t,x))\nn\\&
\ge \theta \lambda e^{-\lambda t}(1+\mid
x\mid^{2\gamma+2})-\frac{1}{2}\theta e^{-\lambda t}
Tr\Big[\sigma.\sigma^\top(t,x)D^{2}_{xx}\mid x\mid^{2\gamma+2}\Big] \nonumber\\
&\qquad \qquad -C(1+|x|)|\nabla_x(-\theta
e^{-\lambda t}\mid x\mid^{2\gamma+2})\sigma(t,x))|\nn\\&
\ge \theta \lambda e^{-\lambda t}(1+\mid
x\mid^{2\gamma+2})-\frac{1}{2}\theta e^{-\lambda t}(C_\sigma)^2(2\gamma+2)(2\gamma+1)|x|^{2\gamma}
\nonumber\\
&\qquad \qquad -C.C_\sigma (1+|x|)(2\gamma+2)\theta
e^{-\lambda t}\mid x\mid^{2\gamma+1}\label{supersol}
\end{align}
where $C_\sigma$ is the constant of boundedness of $\sigma$. But there exists a constant $\lambda_0>0$ such that for any $\lambda \ge \lambda_0$, the right-hand side is positive for any $\theta>0$. Consequently for any $\lambda \ge \lambda_0$ and $\theta>0$, 
$(v+\theta e^{-\lambda t}(1+\mid
x\mid^{2\gamma+2})$ is a supersolution of \eqref{kkk}. \qed
\medskip

In the follwowing lemma, for which we omit the proof since it is classical, we transform the PDE \eqref{kkk} into another one which is more adapted to show uniqueness of the solution of \eqref{kkk}. 
\begin{lemma}\lb{lemmecompar}Let $v(t,x)$ be  an $\R$-valued locally bounded function defined on $\esp$. The function $v$ is a viscosity subsolution (resp. supersolution) of (\ref{kkk}) if and only if
$\bar{v}(t,x)=e^tv(t,x)$, $\txt$, is a viscosity subsolution (resp. supersolution) of the following PDE with obstacle: 
\begin{equation}\label{kkk12}\left\{
\begin{array}{l}
\min\Big\{\bar{v}(t,x)-e^th(t,x),\\
\qquad -{\partial_t \bar{v}}(t,x)+\bar{v}(t,x)-\cl \bar{v}(t,x)-e^t\Phi(t,x,e^{-t}\nabla_x\bar{v}(t,x)\sigma(t,x))\Big\}=0, (t,x)\in [0,T)\times \R^d;\nn\\
\bar{v}(T,x)=e^Tg(x).\qed\nn
\end{array} \right.
\end{equation}
\end{lemma}
 We now address the question of comparison of subsolutions and supersolutions of the PDE \eqref{kkk}. 
\begin{proposition}\label{uni}
Assume that $\Phi$ verifies \bf{(H$\Phi$)}-i), ii). 
Let $\mbox{\underline u}$ (resp. $\underbar v$) be a subsolution (resp. supersolution) of \eqref{kkk}. If $\underbar u,\underbar v$ belong to $\pg$, then $\underbar u\le \underbar v$.\end{proposition}
\nd \textbf{Proof.} First let $\gamma$ be a positive constant such that 
$$|\underline u (t,x)|+|\underline v(t,x)|\le C(1+|x|^{2\gamma}).$$ 
We now that there exists $\lambda_0$ such that for any $\lambda\ge \lambda_0$ and $\theta>0$ such that 
$\underline v(t,x)+\theta e^{-\l t}(1+|x|^{2\g +2})$ still a supersolution of \eqref{kkk}. Therefore it is enough to show that 
$\underline u(t,x)\le \underline v(t,x)+\theta e^{-\l t}(1+|x|^{2\g +2})$ and then to take the limit as $\theta \rw 0$ to obtain the desired result. Next the growth condition on $\underline u$ and $\underline v$ implies the existence of  a positive constant $R$ such that for any $t\in [0,T]$, $|x|\ge R$, $\underline u(t,x)-(\underline v(t,x)+\theta e^{-\l t}(1+|x|^{2\g +2}))<0$. Finally by Lemma \ref{lemmecompar}, $e^t\underline u$ (resp. $e^t(\underline v+\theta e^{-\l t}(1+|x|^{2\g +2}))$) is a viscosity subsolution (resp. supersolution) of \eqref{kkk1}
such that 
$$|e^t\underline u (t,x)|\le C(1+|x|^{2\gamma}) \mbox{ and }|e^t(\underline v(t,x)+\theta e^{-\l t}(1+|x|^{2\g +2}))|\le C(1+|x|^{2\gamma+2}). $$Therefore to obtain the proof it is enough to show that if $u$ (resp. $ w$) is a subsolution (resp. supersolution) of \eqref{kkk1} such that:

i) There exits $R>0$ such that $u(t,x)-w(t,x)<0$ for any $|x| \ge R$ and $t\in [0,T]$; 

ii) $$|u (t,x)|\le C(1+|x|^{2\gamma}) \mbox{ and }|w(t,x)|\le C(1+|x|^{2\gamma+2}). $$

\nd Then $u(t,x)\le w(t,x)$, for any $\txt$.

We will proceed by contradiction and suppose that there exists $(t_0,x_0)\in \spo$ such that 
$u(t_0,x_0)-w(t_0,x_0)>0$ ; w.l.o.g we assume $u$ usc and $w$ lsc. So let $(\bar t,\bar x)$ be such that:
\begin{eqnarray}\label{comp_uni}
 \max_{(t,x) \in [0,T]\times
\mathbb{R}^{d}}(u(t,x)-w(t,x))&=&
\max_{(t,x) \in [0,T[\times
 B_{R}}(u(t,x)-w(t,x))\nonumber\\
&=&(u(\bar{t},\bar{x})-w(\bar{t},\bar{x}))=\eta>0,
\end{eqnarray}
where $B_R := \{x\in \mathbb{R}^{d}; |x|<R\}$ and $(\bar{t},\bar{x}) \in [0,T[\times B_R$.\\

Now let us take
$\btheta$ and $\beta \in (0,1]$, and w.l.o.g we assume $\gamma\ge 2$. Then, for a small
$\epsilon>0$, let us define:

\begin{equation}
\label{phi}
\Phi_{\epsilon}(t,x,y)=u(t,x)-w(t,y)-\frac{1}{2\epsilon}|x-y|^{2\gamma}
-\btheta( |x-\bar{x}|^{2\gamma+2}+|y-\bar{x}|^{2\gamma+2})-\beta
(t-\bar{t})^2.
\end{equation}
Since $u$ is usc and $w$ is lsc, then there exists a
$(t_{\epsilon},x_{\epsilon},y_{\epsilon})\in [0,T]\times \bar B_R
\times \bar B_{R} $ such that:
$$\Phi_{\epsilon}(t_{\epsilon},x_{\epsilon},y_{\epsilon})=\max\limits_{(t,x,y)\in [0,T]\times \bar{B}_{R}\times
\bar{B}_{R}}\Phi_{\epsilon}(t,x,y)$$where $\bar{B}_{R}$ is the closure of $B_R$.
Therefore from the inequality 
$2\Phi_{\epsilon}(t_{\epsilon},x_{\epsilon},y_{\epsilon})\geq
\Phi_{\epsilon}(t_{\epsilon},x_{\epsilon},x_{\epsilon})+\Phi_{\epsilon}(t_{\epsilon},y_{\epsilon},y_{\epsilon})$,
we deduce
\begin{equation}
\frac{1}{\epsilon}|x_{\epsilon} -y_{\epsilon}|^{2\gamma} \leq
(u(t_{\epsilon},x_{\epsilon})-u(t_{\epsilon},y_{\epsilon}))+(w(t_{\epsilon},x_{\epsilon})-w(t_{\epsilon},y_{\epsilon})).
\end{equation}
Consequently $\frac{1}{\epsilon}|x_{\epsilon}
-y_{\epsilon}|^{2\gamma}$ is bounded (thanks to the growth conditions on $u$ and $w$), and as $\epsilon\rightarrow 0$,
$|x_{\epsilon} -y_{\epsilon}|\rightarrow 0$.  By the boundedness of the sequences, one can substract subsequences which we still index by $\epsilon$ such that  $(x_{\epsilon})_\epsilon$ (resp.  $(y_{\epsilon})_\epsilon$, resp. $(t_{\epsilon})_\epsilon$) converges to  $\underbar x$ (resp. $\underbar x$, resp. $\underbar t$) when $\epsilon \rw 0$. Next
 \begin{equation}\label{422}u(\bar{t},\bar{x})-w(\bar{t},\bar{x}) \leq
 \Phi_{\epsilon}(t_{\epsilon},x_{\epsilon},y_{\epsilon})\leq u(t_{\epsilon},x_\epsilon)-w(t_{\epsilon},y_\epsilon).\end{equation}
As $u$ is usc and $w$ is lsc, then we have:
\begin{equation}\begin{array}{l}u(\bar{t},\bar{x})-w(\bar{t},\bar{x}) \leq
 \liminf_{\epsilon \rw 0}\Phi_{\epsilon}(t_{\epsilon},x_{\epsilon},y_{\epsilon})\leq
\limsup_{\epsilon \rw 0}\Phi_{\epsilon}(t_{\epsilon},x_{\epsilon},y_{\epsilon})\leq\\ \qquad\qquad \limsup_{\epsilon \rw 0}(u(t_{\epsilon},x_\epsilon)-w(t_{\epsilon},y_\epsilon))\leq u(\underbar t, \underbar x)-w(\underbar t, \underbar x)\le u(\bar{t},\bar{x})-w(\bar{t},\bar{x}).\end{array}\end{equation}
It follows that: 

i) the sequence $(\Phi_{\epsilon}(t_{\epsilon},x_{\epsilon},y_{\epsilon}))_\eps$ is convergent to $u(\bar{t},\bar{x})-w(\bar{t},\bar{x})$; 

ii) 
$\limsup_{\epsilon \rw 0}(u(t_{\epsilon},x_\epsilon)-w(t_{\epsilon},y_\epsilon))= u(\bar{t},\bar{x})-w(\bar{t},\bar{x})$ and then taking into account of 
\eqref{422}, we deduce that $\lim_{\epsilon \rw 0}(u(t_{\epsilon},x_\epsilon)-w(t_{\epsilon},y_\epsilon))= u(\bar{t},\bar{x})-w(\bar{t},\bar{x})$; 

iii) From \eqref{phi}, we deduce that
 \begin{equation}\label{subsequence}
 (t_\epsilon,x_\epsilon,y_\epsilon)\rightarrow_\eps (\bar{t},\bar{x},\bar{x}).
 \end{equation}
Next $$
\lim_{\epsilon \rw 0}(u(t_{\epsilon},x_\epsilon)-w(t_{\epsilon},y_\epsilon))=u(\bar{t},\bar{x})-w(\bar{t},\bar{x})=\eta>0.$$
Therefore there exists a subsequence of $(\epsilon)$ such that 
$$
(u(t_{\epsilon},x_\epsilon)-w(t_{\epsilon},y_\epsilon))\geq \frac{\eta}{2}.$$
But $w$ is a supersolution then $w(t_\epsilon,y_\epsilon)\geq e^{t_\epsilon}h(t_\epsilon,y_\epsilon)$ and by continuity of $h$ one can find a subsequence such that 
$|e^{t_\epsilon}h(t_\epsilon,y_\epsilon)-e^{t_\epsilon}h(t_\epsilon,x_\epsilon)|<\frac{\eta}{4}$. Therefore for this last subsequence it holds 
\begin{equation}\lb{conditionlimitpouru}u(t_{\epsilon},x_\epsilon)\geq 
e^{t_\epsilon}h(t_\epsilon,x_\epsilon)+\frac{\eta}4.\end{equation}
To proceed let us consider this latter subsequence and let us denote by
\begin{equation}
\varphi_{\epsilon}(t,x,y)=\frac{1}{2\epsilon}|x-y|^{2\gamma}
+\btheta( |x-\bar{x}|^{2\gamma+2}+|y-\bar{x}|^{2\gamma+2})+\beta
(t-\bar{t})^2.
\end{equation}
Then we have: \be
\left\{
\begin{array}{lllll}\label{derive}
D_{t}\varphi_{\epsilon}(t,x,y)=2\beta(t-\bar{t}),\\
D_{x}\varphi_{\epsilon}(t,x,y)= \frac{\gamma}{\epsilon}(x-y)|x-y|^{2\gamma-2} +\btheta({2\gamma}+2)
(x-\bar{x})|x-\bar{x}|^{2\gamma}, \\
D_{y}\varphi_{\epsilon}(t,x,y)= -\frac{\gamma}{\epsilon}(x-y)|x-y|^{{2\gamma}-2} +
\btheta({2\gamma}+2)(y-\bar{x})|y-\bar{x}|^{2\gamma},\\\\
B(t,x,y):=D_{x,y}^{2}\varphi_{\epsilon}(t,x,y)=\frac{1}{\epsilon}
\begin{pmatrix}
a_1(x,y)&-a_1(x,y) \\
-a_1(x,y)&a_1(x,y)
\end{pmatrix}+ \begin{pmatrix}
a_2(x)&0 \\
0&a_2(y)
\end{pmatrix} \\\\
\mbox{with } a_1(x,y)=\gamma|x-y|^{{2\gamma}-2}I+{\gamma}({2\gamma}-2)(x-y)(x-y)^*|x-y|^{{2\gamma}-4}  \mbox{ and }\\
a_2(x)=\btheta({2\gamma}+2)|x-\bar{x}|^{2\gamma}I+2\btheta\gamma({2\gamma}+2)
(x-\bar{x})(x-\bar{x})^*|x-\bar{x}|^{2\gamma-2} .
\end{array}
\right. \ee Applying now 
the result by Crandall et al. (Theorem 8.3, {\cite{CIL}) to the
function $$
u(t,x)-w(t,y)-\varphi_{\epsilon}(t,x,y) $$ at
the point $(t_\epsilon,x_\epsilon,y_\epsilon)$ (we choose $\epsilon$ small enough in such a way that $t_\eps <T$, and $|x_\eps|<R$ and $|y_\eps|<R$), for any $\upsilon>0$ we can find
$c,c_1 \in \R$, $q_1,q_2 \in \R^d$ and $X,Y \in S_d$, such that: 

\be \label{lemmeishii}
\left\{
\begin{array}{lllll}
(c,q_1,X)
\in \bar J^{2,+}u(t_\epsilon,x_\epsilon) \mbox{ and }
(-c_1,q_2,Y)\in \bar J^{2,-}
w(t_\epsilon,y_\epsilon),\\
q_1=D_{x}\varphi_{\epsilon}(t_\epsilon,x_\epsilon,y_\epsilon)=
\frac{{\gamma}}{\epsilon}(x_\epsilon-y_\epsilon)|x_\epsilon-y_\epsilon|^{2\gamma-2} +\theta({2\gamma}+2)
(x_\epsilon-\bar{x})|x_\epsilon-\bar{x}|^{2\gamma},\\
q_2= -D_{y}\varphi_{\epsilon}(t_\epsilon,x_\epsilon,y_\epsilon)=\frac{{\gamma}}{\epsilon}(x_\epsilon-y_\epsilon)|x_\epsilon-y_\epsilon|^{2\gamma-2} -
\theta({2\gamma}+2)(y_\epsilon-\bar{x})|y_\epsilon-\bar{x}|^{{2\gamma}},\\
c+c_1=D_{t}\varphi_{\epsilon}(t_\epsilon,x_\epsilon,y_\epsilon)=2\beta(t_\epsilon-\bar{t}) \mbox{ and finally }\\
-(\frac{1}{\upsilon}+||B(t_\epsilon,x_\epsilon,y_\epsilon)||)I\leq
\begin{pmatrix}
X&0 \\
0&-Y
\end{pmatrix}\leq B(t_\epsilon,x_\epsilon,y_\epsilon)+\upsilon B(t_\epsilon,x_\epsilon,y_\epsilon)^2.
\end{array}
\right. \ee Taking now into account (\ref{conditionlimitpouru}), and the
definition of viscosity solution, we get:
$$\begin{array}{l}-c-\frac{1}{2}Tr[\sigma^\top(t_\epsilon,x_\epsilon)X\sigma(t_\epsilon,x_\epsilon)]
+u(t_\epsilon,x_\epsilon)-e^{t_\epsilon}\Phi(t_\epsilon,x_\epsilon,e^{-t_\epsilon}q_1\sigma(t_\epsilon,x_\epsilon))\leq 0\\
\mbox{ and
}\\c_1-\frac{1}{2}Tr[\sigma^\top(t_\epsilon,y_\epsilon)Y\sigma(t_\epsilon,y_\epsilon)]
+ w(t_\epsilon,y_\epsilon)-e^{t_\epsilon}\Phi(t_\epsilon,y_\epsilon,e^{-t_\epsilon}q_2\sigma(t_\epsilon,y_\epsilon))\geq
0.\end{array}$$  Then
\begin{equation}
\begin{array}{llll}
\label{viscder}
u(t_\epsilon,x_\epsilon)-w(t_\epsilon,y_\epsilon)-c-c_1&\leq \frac{1}{2}Tr[\sigma^\top(t_\epsilon,x_\epsilon)X\sigma(t_\epsilon,x_\epsilon)-\sigma^\top(t_\epsilon,y_\epsilon)Y\sigma(t_\epsilon,y_\epsilon)]

\\&\qquad+e^{t_\epsilon}\Phi(t_\epsilon,x_\epsilon,e^{-t_\epsilon}q_1\sigma(t_\epsilon,x_\epsilon))-e^{t_\epsilon}\Phi(t_\epsilon,y_\epsilon,e^{-t_\epsilon}q_2\sigma(t_\epsilon,y_\epsilon)).
\end{array}
\end{equation}
But from (\ref{derive}) there exist two constants $C$ and $C_1$ (which may change from line to line) such
that:
$$||a_1(x_\epsilon,y_\epsilon)||\leq C|x_\epsilon- y_\epsilon|^{2\gamma-2} \mbox{ and }(||a_2(x_\epsilon)||\vee ||a_2(y_\epsilon)||)\leq C_1 \btheta.$$
On the other hand we have $(B:= B(t_\epsilon,x_\epsilon,y_\epsilon))$
$$B\leq \frac{C}{\epsilon}|x_\epsilon - y_\epsilon|^{2\g-2}
\begin{pmatrix}
I&-I \\
-I&I
\end{pmatrix}+ C_1 \theta I.$$
It follows that:
\begin{equation}
B+\upsilon B^2 \leq C(\frac{1}{\epsilon}|x_\epsilon - y_\epsilon|^{2\gamma-2}+
\frac{\epsilon_1}{\epsilon^2}|x_\epsilon - y_\epsilon|^{4\gamma-4})\begin{pmatrix}
I&-I \\
-I&I
\end{pmatrix}+ C_1\btheta I.
\end{equation}
Choosing now $\upsilon=\epsilon$, yields:
\begin{equation}
\label{ineg_matreciel}
B+\eps B^2 \leq \frac{C}{\epsilon}(|x_\epsilon - y_\epsilon|^{2\gamma-2}+|x_\epsilon - y_\epsilon|^{4\gamma-4})\begin{pmatrix}
I&-I \\
-I&I
\end{pmatrix}+ C_1\btheta I.
\end{equation}
From the Lipschitz continuity of $\sigma$, (\ref{lemmeishii}) and
(\ref{ineg_matreciel}) we have:
\begin{equation}\label{termensigma}\frac{1}{2}Tr[\sigma^\top(t_\epsilon,x_\epsilon)X\sigma(t_\epsilon,x_\epsilon)-\sigma^\top(t_\epsilon,y_\epsilon)
Y\sigma(t_\epsilon,y_\epsilon)]\leq \frac{C}{\epsilon}(|x_\epsilon - y_\epsilon|^{2\gamma}+|x_\epsilon - y_\epsilon|^{4\gamma-2}) +C_1 \btheta.\end{equation} Next taking into account of \bf{(H$\Phi$)} we have:
\begin{align}
&e^{t_\epsilon}\Phi(t_\epsilon,x_\epsilon,e^{-t_\epsilon}q_1\sigma(t_\epsilon,x_\epsilon))-e^{t_\epsilon}\Phi(t_\epsilon,y_\epsilon,e^{-t_\epsilon}q_2\sigma(t_\epsilon,y_\epsilon))\nn\\
&\le e^{t_\epsilon}\{C_R(|e^{-t_\epsilon}q_1|+
|e^{-t_\epsilon}q_2|)|x_\epsilon-y_\epsilon|+
C_R|e^{-t_\epsilon}q_1-e^{-t_\epsilon}q_2|+\Psi_R(|x_\epsilon-y_\epsilon|)\}\nn
\end{align}
But $\lim_{\eps \rw 0}(|q_1|+|q_2|)|x_\epsilon-y_\epsilon|=0$ and 
$\lim_{\eps \rw 0}|q_1-q_2|=0$.
It follows that 
$$
\limsup_{\eps \rw 0}e^{t_\epsilon}\Phi(t_\epsilon,x_\epsilon,e^{-t_\epsilon}q_1\sigma(t_\epsilon,x_\epsilon))-e^{t_\epsilon}\Phi(t_\epsilon,y_\epsilon,e^{-t_\epsilon}q_2\sigma(t_\epsilon,y_\epsilon))\le 0.
$$
Next go back to \eqref{viscder} take the superior limit w.r.t $\eps$ on each hand-side and take into account of \eqref{termensigma} to obtain: 
$$
u(\bar{t},\bar{x})-w(\bar{t},\bar{x})\le C_1\btheta.
$$
Send now $\btheta
\rightarrow0$ to obtain that $u(\bar{t},\bar{x})-w(\bar{t},\bar{x})=0$ which is contradictory. The proof is now complete. \qed

As a by-product we have

\begin{corollary}Under \bf{(H$\Phi$)}, if the PDE \eqref{kkk} has a solution in $\pg$, then it is unique and continuous.  \qed
\end{corollary}
Let us now go back to the HJB equation associated with the mixed control problem, i.e., equation \eqref{kkk} when $\Phi$ is replaced with $H^*$ whihc reads 
\begin{equation}\label{hjb}
\left\{
\begin{array}{l}
\min\left[v(t,x)-h(t,x),
-\frac{\partial v}{\partial t}(t,x)-\cl v(t,x)-H^*(t,x,\nabla_xv(t,x)\sigma(t,x))\right]=0, \,\,(t,x)\in[0,T[ \times \R ^d;\\\\
v(T,x)=g(x),\, x\in\R^d.
\end{array}
\right.
\end{equation}
To begin with we will focus on the properties of the function $H^*$.
\begin{lemma}\lb{conth*}${}$\\
i) The function $H^*$ is continuous in $(t,x,z)$. \\
ii) $H^*$ verifies \eqref{cd1gfi} and \eqref{cd2gfi}.
\end{lemma}
\noindent \bf{Proof}: i) Let $(t,x,z)$ and $(t',x',z')$ be fixed. Without loss of generality we assume that 
$|t-t'|+|x-x'|+|z-z'|\le 1$. An easy computation shows that
\begin{align}|H^*(t,x,z)-H^*(t',x',z')|&=|H^*(t,x,z)-H^*(t',x',z)+H^*(t',x',z)-H^*(t',x',z')|\nn\\
&\le |z|\sup_{a\in \ca}|\sigma^{-1}(t,x)f(t,x,a)-\sigma^{-1}(t',x')f(t',x',a)|\nn\\&\nn\qquad \qquad\qquad+\sup_{a\in \ca}|\Gamma(t,x,a)-\Gamma(t',x',a)| +C(1+|x'|)|z-z'|.
\end{align}
Next as $f$, $\sigma^{-1}$ and $\G$ are continuous and $\ca$ is compact then the right-hand side of the previous inequality goes to $0$ when $(t',x',z')\rw (t,x,z)$. Thus $H^*$ is continuous. 

\ni ii) \begin{align}\nn|H^*(t,x,z)-H^*(t,x,z')|&=|\sup_{a\in \ca}\{z\sigma^{-1}(t,x)f(t,x,a)+\Gamma(t,x,a)\}-\sup_{a\in \ca}\{z'\sigma^{-1}(t,x)f(t,x,a)+\Gamma(t,x,a)\}|\\&\nn\le
\sup_{a\in \ca}|z-z'||\sigma^{-1}(t,x)f(t,x,a)|\\&\le C(1+|x|)|z-z'|\nn
\end{align}
since $\sigma^{-1}$ is bounded and $f$ of linear growth. Thus $H^*$ verifies \eqref{cd1gfi}. Finally let us show that $H^*$ verifies \eqref{cd2gfi}. Let $\kappa$ be fixed and $(t,z,z')\in [0,T]\times \R^{2d}$, $x,y\in \R^d$ such that $|x|\le \kappa$ and $|y|\le \kappa$.
\begin{align}&H^*(t,x,z\sigma(t,x))-H^*(t,y,z'\sigma(t,y))\nn\\&\qquad \qquad =\sup_{a\in \ca}\{zf(t,x,a)+\Gamma(t,x,a)\}-\sup_{a\in \ca}\{z'f(t,y,a)+\Gamma(t,y,a)\}\nn\\&\qquad \qquad \le \sup_{a\in \ca}\{zf(t,x,a)-z'f(t,y,a)\}+\sup_{a\in \ca}\{\Gamma(t,x,a)-\Gamma(t,y,a)\}\nn\\&\qquad \qquad \le C(1+|x|)|z-z'|+\sup_{a\in \ca}\{z'(f(t,x,a)-f(t,y,a))\}+\sup_{a\in \ca}\{\Gamma(t,x,a)-\Gamma(t,y,a)\}
\nn\\&\qquad \qquad \le C(1+|x|)|z-z'|+C|z'||x-y|+
\Psi^\G_\kappa(|x-y|). 
\end{align}
where $\Psi^\G_\kappa$ is the modulus of continuity 
of $\G$ on $[0,T]\times \bar B(0,\kappa)\times \ca$ 
($\bar B(0,\kappa)$ is the closure in $\R^d$ of the open ball centered in 0 and of radius $\kappa$). In the last inequality we have used the fact that $f$ is Lipschitz w.r.t $x$. The proof now follows since $|x|\le \kappa$.\qed
\medskip

We now have the following result related to $\bar u^m$.
\begin{proposition}\lb{contbarum}For any $m\ge 0$, $\bar u^m$ is continuous and is the unique viscosity solution in $\Pi_g$ of the following PDE with obstacle:
\begin{equation}\label{vis-m}\left\{\begin{array}{l}\min\left[\bar{u}_{m}(t,x)-h(t,x),\right.\\\\\left.-\partial_t\bar{u}_{m}(t,x)-\cl\bar{u}_{m}(t,x)-\bar{H}^{*m}(t,x,\sigma(t,x)\nabla_x\bar{u}_{m}(t,x))\right]=0,(t,x)\in[0,T)\times \R^d;\\\\
\bar u(T,x)=g(x).\end{array}\right.\end{equation}
\end{proposition}
\ni\bf{Proof}: First recall tha by \eqref{croipolyumu}, $\bar u^m$ belongs to $\pg$. Next let us show that $\bar{u}_m$ is a supersolution  of (\ref{vis-m}). The function $\bar u^m$ is lsc and $\bar u^m(T,x)=g(x)$. Thus the terminal condition is verified. On the other hand, $\bar u^m_*=\bar u^m$. So let $t<T$ and $(p,q,X)\in \bar J^{2,-} \bar{u}_m(t,x)$. As $\bar u^m=\lim_n\nearrow \bar u_{n,m}$ then thanks to Lemma 6.1 in \cite{CIL}, there exist sequences:
$$
\begin{array}{l}
n_j \rightarrow+\infty,\quad (t_j,x_j)\rightarrow (t,x), \quad (p_j,q_j,X_j)\in
\bar J^{2,-}\bar{u}_{n_j,m}(t_j,x_j),
\end{array}
$$
such that $$(p_j,q_j,X_j)\rightarrow (p,q,X).$$ But for any $j$,
$$
\begin{array}{l}
-p_j -\frac{1}{2}Tr[\sigma^\top(t_j,x_j) X_j \sigma(t_j,x_j)] \ge \bar{H}^{*n_j,m}(t_j,x_j,q_j\sigma(t_j,x_j))\end{array}
$$
since $\bar u_{n,m}$ is a viscosity solution of \eqref{edpunm}. But by Dini's Theorem,  
$(\bar{H}^{*n,m})_n$ converges uniformly to $\bar{H}^{*m}$ on compact subsets. Therefore take the limit w.r.t. $j$ in each hand-side of the previous inequality to obtain 
$$
\begin{array}{l}
-p -\frac{1}{2}Tr[\sigma^\top(t,x) X \sigma(t,x)] \ge \bar{H}^{*m}(t,x,q\sigma(t,x)).\end{array}
$$
Finally as $\bar u^m(t,x)=\bar Y_t^{m,t,x}\ge h(t,x)$, then $\bar u^m$ is a viscosity supersolution of \eqref{vis-m}.

Let us now show that $\bar u^{m*}$ is a subsolution of \eqref{vis-m}. Let $(t,x)\in [0,T)\times \R^d$. We obvioulsy have $\bar {u}^{m*}(t,x)\ge h(t,x)$ since $h$ is continuous. So assume that $\bar {u}^{m*}(t,x)> h(t,x)$ and let 
$(p,q,X)\in J^{2,-}\bar{u}^{m*}(t,x).$ As 
$\bar u^m=\lim_n\nearrow \bar u_{n,m}$ and $\bar u_{n,m}$ is continuous then 
$$\bar u^{m*} = \lim_{n\rw \infty}{\sup}^*\bar u_{n,m}$$ 
where 
$$
\lim_{n\rw \infty}{\sup}^*\bar u_{n,m}(t,x)=\limsup_{n\rw \infty,(t',x')\rw (t,x), t'<T}\bar u_{n,m}(t',x')
$$
(see \cite{barles}, pp.91).
Then once more by Lemma 6.1 in \cite{CIL}, there exist sequences such that
$$
\begin{array}{l}
n_j \rightarrow+\infty,\quad (t_j,x_j,\bar{u}_{n_j,m}(t_j,x_j))\rightarrow (t,x,\bar{u}^{m*}(t,x)),\quad (p_j,q_j,X_j)\in
J^{2,-}\bar{u}_{n_j,m}(t_j,x_j),
\end{array}
$$
and $$(p_j,q_j,X_j)\rightarrow (p,q,X).$$ But there exists a subsequence which we still index by $j$ such that for any $j$,
$\bar{u}_{n_j,m}(t_j,x_j)>h(t_j,x_j)$ since $\bar {u}^{m*}(t,x)> h(t,x)$. The subsolution property of $\bar{u}^{n_j,m}$ implies that:
$$
\begin{array}{l}
-p_j -\frac{1}{2}Tr[\sigma^\top (t_j,x_j)X_j \sigma(t_j,x_j)] \le \bar{H}^{*n_j,m}(t_j,x_j,q_j\sigma(t_j,x_j)). 
\end{array}
$$
Hence, as for the supersolution property, taking the limit as $j\rightarrow +\infty$, we conclude that:
$$-p -\frac{1}{2}Tr[\sigma^\top (t,x)X \sigma(t,x)] -\bar{H}^{*m}(t,x,q\sigma(t,x))\leq 0, $$
i.e. $\bar{u}_m$ verifies the subsolution property in $(t,x)\in [0,T)\times \R^d$. It remains to show the terminal condition is satisfied i.e.
$\bar u^{m*}(T,x)=g(x)$. It is classical, however we give it in its main steps for completeness. First note that for any $x\in\R^d$,  \be \label{plusgrand}\bar u^{m*}(T,x)\geq
\bar u_{n_0,m}(T,x)=g(x).\ee Assume now that for some $x_0\in \R^d$, 
 \be \label{plusgrand2}\bar u^{m*}(T,x_0)-g(x_0)=2\eps>0\ee
and let us
construct a contradiction. Let $(t_k,x_k)_{k\geq 1}$ be a sequence
in $[0,T)\times \R^d$ such that:
$$(t_k,x_k)\rightarrow (T,x_0) \mbox{ and
}\bar u^{m}(t_k,x_k)\rightarrow \bar u^{m*}(T,x_0) \mbox{ as }k\rightarrow
\infty.
$$
Since $\bar u^{m*}$ is $usc$ and of polynomial growth and taking into
account of (\ref{plusgrand}) and the inequalilty $g(x_0)\ge h(T,x_0)$, we can find a sequence $(\varrho^n)_{n\geq
0}$ of functions of $\cC^{1,2}(\esp)$ such that $\vro^n\rightarrow
\bar u^{m*}$ and, on some neighbourhood $B_n$ of $(T,x_0)$ we have:
\be\label{eq7x}\min\{\vro^n(t,x)-g(x), \vro^n(t,x)-h(t,x)\}\geq \eps,\,\,\forall (x,t)\in
B_n.\ee After possibly passing to a subsequence of $(t_k,x_k)_{k\geq
1}$ we can then assume that it holds on $B^n_k:=[t_k,T]\times
B(x_k,\delta_n^k)$ for some $\delta^n_k \in (0,1)$ small enough in such a
way that $B^n_k\subset B$. Now since $\bar u^{m*}$ is locally bounded
then there exists $\zeta >0$ such that $|\bar u^{m*}|\leq \zeta$ on
$B_n$. We can then assume that $\vro ^n\geq -2\zeta$ on $B_n$. Next
let us define:$$ \tilde \vro^n_k(t,x):=\vro^n(t,x)+\frac{4\zeta
|x-x_k|^2}{(\delta^n_k)^2}+\sqrt{T-t}.$$ Note that $\tilde \vro^n_k\geq
\vro^n$ and \be\label{eq8}(\bar u^{m*}-\tilde \vro^n_k)(t,x)\leq -\zeta
\mbox{ for }(t,x)\in [t_k,T]\times
\partial B(x_k,\delta^n_k).\ee
Next since $\partial_t(\sqrt{T-t})\rightarrow -\infty$ as
$t\rightarrow T$, we can choose $t_k$ large enough in front of
$\delta^n_k$ and the derivatives of $\vro^n$ to ensure that
\be\label{eq9}-\cl \tilde \vro^n_k (t,x)\geq 0 \mbox{ on }B_n^k.\ee
Next let us consider the following stopping time
$\theta_n^k:=\inf\{s\geq t_k, (s,X^{t_k,x_k}_s)\in {B_n^k}^c\}\wedge
T$ where ${B_n^k}^c$ is the complement of ${B_n^k}$, and
$\vartheta_k:=\inf\{s\geq t_k, \bar u^{m*}(s,X^{t_k,x_k}_s)=h(s,X^{t_k,x_k}_s)\}\wedge T.$
Using It\^o's formula and taking into account (\ref{eq7x}),
(\ref{eq8}) and (\ref{eq9}) to obtain:
$$\ba{ll}
\tilde \vro^n_k(t_k,x_k)&=\E[\tilde \vro^n_k(\theta_n^k\wedge
\vartheta_k,X^{t_k,x_k}_{\theta_n^k \wedge
\vartheta_k})-\int_{t_k}^{\theta_n^k\wedge \vartheta_k}
\cl \tilde \vro^n_k (r,X_r^{t_k,x_k})dr]\\
{}& \geq \E[\tilde
\vro^n_k(\theta_n^k,X^{t_k,x_k}_{\theta_n^k})\ind_{[\theta_n^k\leq
\vartheta_k]}
+\tilde\vro^n_k(\vartheta_k,X^{t_k,x_k}_{\vartheta_k})\ind_{[\vartheta_k<
\theta_n^k]}]\\{}&= \E[\{\tilde
\vro^n_k(\theta_n^k,X^{t_k,x_k}_{\theta_n^k})\ind_{[\theta_n^k<T]}+
\tilde
\vro^n_k(T,X^{t_k,x_k}_{T})\ind_{[\theta_n^k=T]}\}\ind_{[\theta_n^k\leq
\vartheta_k]}
+\tilde\vro^n_k(\vartheta_k,X^{t_k,x_k}_{\vartheta_k})\ind_{[\vartheta_k<
\theta_n^k]}]
\\
{}&\geq
\E[\{(\bar u^{m*}(\theta_n^k,X^{t_k,x_k}_{\theta_n^k})+\zeta)\ind_{[\theta_n^k<T]}+
(\eps
+g(X^{t_k,x_k}_{T}))\ind_{[\theta_n^k=T]}\}\ind_{[\theta_n^k\leq
\vartheta_k]} \\{}&\qquad \qquad +\{\eps+h(\vartheta_k,X^{t_k,x_k}_{\vartheta_k})\}\ind_{[\vartheta_k<
\theta_n^k]}]\\&\geq
\E[\{(\bar u^{m}(\theta_n^k,X^{t_k,x_k}_{\theta_n^k})+\zeta)\ind_{[\theta_n^k<T]}+
(\eps
+g(X^{t_k,x_k}_{T}))\ind_{[\theta_n^k=T]}\}\ind_{[\theta_n^k\leq
\vartheta_k]} \\{}&\qquad \qquad +\{\eps+h(\vartheta_k,X^{t_k,x_k}_{\vartheta_k})\}\ind_{[\vartheta_k<
\theta_n^k]}]
\\
{}&\geq \E[\bar u^{m}(\theta_n^k\wedge
\vartheta_k,X^{t_k,x_k}_{\theta_n^k \wedge
\vartheta_k})]+\zeta\wedge \eps\\{}&
=\E[\bar u^{m}(t_k,x_k)-\int_{t_k}^{\theta_n^k\wedge
\vartheta_k}\bar H^{*m}(s,X_s^{t_k,x_k},\bar Z^{t_k,x_k,m}_s)ds]+\zeta\wedge
\eps \ea$$ since on $[t_k,\vartheta_k]$, $d\bar K^{t,x,m}=0$. Finally
since $\bar u^{m}$ and $\|\bar Z^{t_k,x_k,m}\|_{{\cal H}^{2,d}}$ are of polynomial growth (Lemma \ref{lemmecroissancepoly}) we deduce that
$\lim_{k\rightarrow \infty}\E[\int_{t_k}^{\theta_n^k\wedge
\vartheta_k}\bar H^{*m}(s,X_s^{t_k,x_k},\bar Z^{t_k,x_k,m}_s)ds]=0$. Therefore taking the limit in the previous
inequalities yields:
$$\label{eq10} \lim_{k\rightarrow
\infty}\tilde \vro^n_k(t_k,x_k)=\lim_{k\rightarrow
\infty}\vro^n(t_k,x_k)+\sqrt{T-t_k}=\vro^n(T,x_0)\geq
\lim_{k\rightarrow \infty} \bar u^m(t_k,x_k)+\zeta\wedge
\eps=\bar u^{m*}(T,x_0)+\zeta\wedge \eps.$$ But this is contradictory
since $\vro^n\rightarrow \bar u^{m*}$ pointwise as $n\rightarrow
\infty$. Thus for any $x\in \R^d$ we have
$\bar u^{m*}(T,x)=g(x)$ and this completes the proof.

Finally to show continuity of $\bar u^m$ and uniqueness of the solution it is enough to show that $\bar H^{*m}$ satisfies the assumptions \bf{(H$\Phi)$}. But for any $(t,x,z,z')$
$$
|\bar H^{*m}(t,x,z)-\bar H^{*m}(t,x,z')|\le 2
|H^*(t,x,z)-H^*(t,x,z')|\le 
2C(1+|x|)|z-z'|$$
since $f$ is of linear growth, $\sigma^{-1}$ is bounded and $|\rho_m|\le 1$. Thus \eqref{cd1gfi} is satisfied. 

On the other hand, 
$$\ba{ll}
\bar H^{*m}(t,x,z\sigma(t,x))-\bar H^{*m}(t,y,z'\sigma(t,y))=\\
\qquad\qquad \underbrace{H^{*+}(t,x,z\sigma(t,x))-H^{*+}(t,y,z'\sigma(t,y))}_{\mbox{(I)}}+\\\qquad\qquad 
\{\underbrace{-(H^{*-}(t,x,z\sigma(t,x))-H^{*-}(t,y,z'\sigma(t,y)))}_{\mbox{(II)}}\rho_m(x)-
\underbrace{H^{*-}(t,y,z'\sigma(t,y))}_{\mbox{(III)}}(\rho_m(x)-\rho_m(y)).
\end{array}$$
Now to conclude it is enough to remark that: i) $H^*$ verifies \eqref{cd2gfi} and to use the inequality 
$u^+-v^+\leq (u-v)^+$ with (I); ii) $H^*$ verifies \eqref{cd2gfi} and to use the inequality 
$u^--v^-\leq (v-u)^+$ with (II); iii) $H^{*-}$ verifies \eqref{cd1gfi} and $\rho_m$ is continuous Lipschitz. 
The proof is now complete.  \qed
\medskip

We are now ready to state the main result of this section.

\begin{theorem}The function $u$ is continuous and is the unique viscosity solution in $\pg$ of the following PDE with obstacle:
\begin{equation}\label{vis-nm}\left\{\begin{array}{l}\min\left[{u}(t,x)-h(t,x),\right.\\\\\left.-\partial_tu(t,x)-\cl u(t,x)-{H}^*(t,x,\sigma(t,x)\nabla_x u(t,x))\right]=0,(t,x)\in[0,T)\times \R^d;\\\\
u(T,x)=g(x).\end{array}\right.\end{equation}
\end{theorem}
\ni\bf{Proof}:  The proof is obtained in the same way as we did for $\bar u^m$ in the previous proposition since: i)  by \eqref{croipolyumu} we now  that $u$ belongs to $\pg$; ii) $\bar u^m$ verifies the PDE \eqref{vis-nm}, is continuous and the sequence 
$(\bar u^m)_m$ is decreasing and converges to $u$; iii) $H^*$ verifies \eqref{cd1gfi} and \eqref{cd2gfi}. The details are left to the care of the reader. \qed

\begin{remark}The characterization of $u$ (see Remark \ref{remarque31}) as the value function of the mixed control problem allows to show directly that $u$ is continuous even if this proof is rather tedious. However there is no way to show that $\bar u^m$ is continuous without using the PDEs as we did previously in Proposition \ref{contbarum}.  \qed
\end{remark}

\end{document}